\documentclass[english,12pt]{smfart}

\usepackage[english]{babel}
\usepackage{amssymb}
\usepackage{xypic}

\def\Z{{\mathbb Z}}
\def\Q{{\mathbb Q}}
\def\Qp{\Q_p}
\def\Zp{\Z_p}
\def\Fp{{\mathbb F}_p}
\def\Qpbar{\overline{\Qp}}
\def\Fpbar{\overline{\Fp}}
\def\O{\mathcal{O}_L}
\def\g{{\rm Gal}(\Qpbar/L)}
\def\gL{{\rm Gal}(L'/L)}
\def\v{{\rm val}_L}
\def\abs{|\ |}

\newtheorem{definit}{\bf{ Definition }}[section]
\newtheorem{prop}[definit]{\bf{ Proposition  }}
\newtheorem{lem}[definit]{\bf{ Lemma  }}
\newtheorem{thm}[definit]{\bf{ Theorem  }}
\newtheorem{cor}[definit]{\bf{ Corollary  }}
\newtheorem{conj}[definit]{\bf{ Conjecture  }}
\newtheorem{rem}[definit]{\bf{ Remark  }}

\author[C. Breuil]{Christophe Breuil}
\address{C.N.R.S. \& I.H.\'E.S.\\
Le Bois-Marie\\
35 route de Chartres\\
91440 Bures-sur-Yvette\\
France} \email{breuil@ihes.fr}

\author[P. Schneider]{Peter Schneider}
\address{Mathematisches Institut\\
Universit\"at M\"unster\\
Einsteinstrasse 62\\
48149, M\"unster\\
Deutschland} \email{pschnei@math.uni-muenster.de}

\title[First steps towards $p$-adic Langlands functoriality]{First steps towards $p$-adic Langlands functoriality}

\begin{document}

\begin{abstract}
By the theory of Colmez and Fontaine, a de Rham representation of the Galois group of a local field roughly corresponds to a representation of the Weil-Deligne group equipped with an admissible filtration on the underlying vector space. Using a modification of the
classical local Langlands correspondence, we associate with any pair consisting of a Weil-Deligne group representation and a type of a filtration (admissible or not) a specific locally algebraic representation of a general linear group. We advertise the conjecture that this pair comes from a de Rham representation if and only if the corresponding locally algebraic representation carries an invariant norm. In the crystalline case, the Weil-Deligne group representation is unramified and the associated locally algebraic representation can be studied using the classical Satake isomorphism. By extending the
latter to a specific norm completion of the Hecke algebra, we show that the existence of an invariant norm implies that our pair, indeed, comes from a crystalline representation. We also show, by using the formalism of Tannakian categories, that this latter fact is compatible with classical unramified Langlands functoriality and therefore generalizes to arbitrary split reductive groups.
\end{abstract}

\begin{altabstract}
Par la th\'eorie de Colmez et Fontaine, une repr\'esentation de de Rham du groupe de Galois d'un corps local correspond essentiellement \`a une repr\'esentation du groupe de Weil-Deligne dont l'espace sous-jacent est muni d'une filtration admissible. En modifiant la correspondance locale de Langlands, on associe \`a chaque couple form\'e d'une repr\'esentation du groupe de Weil-Deligne et des poids d'une filtration (admissible ou pas) une repr\'esentation localement alg\'ebrique particuli\`ere d'un groupe lin\'eaire g\'en\'eral. On conjecture qu'un couple provient d'une repr\'esentation de de Rham si et seulement si la repr\'esentation localement alg\'ebrique correspondante poss\`ede une norme invariante. Dans le cas cristallin, la repr\'esentation du groupe de Weil-Deligne est non-ramifi\'ee et la repr\'esentation localement alg\'ebrique associ\'ee peut s'\'etudier gr\^ace \`a l'isomorphisme de Satake classique. En prolongeant ce dernier \`a une compl\'etion de l'alg\`ebre de Hecke, on montre que l'existence d'une norme invariante comme ci-dessus implique que le couple provient effectivement d'une repr\'esentation cristalline. On montre aussi, en utilisant le formalisme des cat\'egories tannakiennes, que ce dernier fait est compatible avec la fonctorialit\'e de Langlands non-ramifi\'ee classique, et donc qu'il se g\'en\'eralise \`a tout groupe r\'eductif d\'eploy\'e.
\end{altabstract}

\maketitle

\setcounter{tocdepth}{2}

\tableofcontents

\section{Introduction and notations}

The authors strongly believe in the existence of a $p$-adic
extension of the local Langlands correspondence and even the
functoriality principle. Unfortunately, at the present time, there
is not even a precise picture yet which two sets this extension will
put into correspondence. On the Galois side we should at least have
all $p$-adic Weil group representations. On the reductive group side
the evidence points into the direction of a set related to the set
of isomorphism classes of all topologically irreducible admissible
Banach space representations. But it seems very difficult to
construct such Banach space representations. There has been progress
recently by Breuil/Berger and by Colmez only for the group
${\rm GL}_2(\mathbb{Q}_p)$.

If we restrict attention on the Galois side to $(d+1)$-dimensional
de Rham representations of the Galois group $\g$ of some finite
extension $L$ of $\mathbb{Q}_p$ then the first author some time ago
has put forward the following philosophy. By the theory of Colmez
and Fontaine, a de Rham representation gives rise (roughly) to a
vector space $D$ which carries a filtration as well as an action of
the Weil-Deligne group of $L$, the two being in a numerical relation
called admissibility. The classical local Langlands correspondence
associates with the semisimplification of the Weil-Deligne group
action an irreducible smooth representation of ${\rm GL}_{d+1}(L)$. On the
other hand, the type of the filtration can be viewed as a highest
weight giving rise to an irreducible rational representation of
${\rm GL}_{d+1}(L)$. The tensor product of the two forms a locally
algebraic representation $\Pi(D)$. Note that the construction of
$\Pi(D)$ does not make use of the admissibility relation on $D$.
Rather the admissibility of the filtration on $D$ should be
reflected by the existence of a ${\rm GL}_{d+1}(L)$-invariant norm on
$\Pi(D)$. It is one of the purposes of this paper to turn this
philosophy into a precise conjecture which is done in section 4. The
problem, of course, is that the de Rham representation can be such
that the smooth part of $\Pi(D)$ is trivial whereas the algebraic
part is not. Then $\Pi(D)$ cannot carry an invariant norm. In order
to overcome this difficulty we suggest to use a modified version of
the classical local Langlands correspondence. By the Langlands
classification, every irreducible smooth representation is the unique
irreducible quotient of a particular parabolically induced
representation. We propose to use (a twist of) the latter instead of
its irreducible quotient. In this way the smooth part of $\Pi(D)$
always is infinite dimensional and, in fact, reducible in general.
In section 5, we establish some special cases of our conjecture as
well as some partial results in its direction.

If we restrict attention further to crystalline Galois
representations then the smooth part of $\Pi(D)$ is unramified. The
classical local Langlands correspondence for unramified
representations is encapsulated in the Satake isomorphism which
computes the spectrum of the Satake-Hecke algebra of ${\rm
GL}_{d+1}(L)$. The Satake-Hecke algebra together with its universal
module has a natural norm completion leading to a Banach-Hecke
algebra acting on a universal Banach module. This was shown by the
second author together with Teitelbaum in \cite{ST}. It was also
shown in that paper that this Banach-Hecke algebra naturally is the
algebra of analytic functions on an explicitly given affinoid
domain. Moreover, the defining conditions for this affinoid domain
turn out to be equivalent to the admissibility condition for
filtrations. This means that any crystalline Galois representation
satisfying a certain regularity condition on its Hodge-Tate weights
gives rise to a point in such an affinoid domain and hence to a
specialization in that point of the corresponding universal Banach
module. Our earlier conjecture means in this picture that this
Banach space representation of ${\rm GL}_{d+1}(L)$ obtained by
specialization is non-zero. Unfortunately, in \cite{ST} an embedding
of $L$ into the coefficient field $K$ of our representations was
distinguished. This had the consequence that all of the above could
only be shown for a subclass of the crystalline Galois
representations called special ones. Even worse, the corresponding
normalizations in \cite{ST} are, as we believe now, misleading. In
sections 2 and 3, we take up this theory again in a completely
general way. By working systematically with $\mathbb{Q}_p$-rational
representations of ${\rm GL}_{d+1}(L)$, we do obtain the above
results for arbitrary ``regular'' crystalline Galois
representations.

At this point it should be stressed that, in this paper, we do not set
up an actual conjectural correspondence between Galois and Banach
space representations. In both pictures, our conjecture is of the
form that on a certain object there exists an admissible filtration
if and only if a certain Banach space representation associated to
this object (either by completion or by specialization) is non-zero.
But with the exception of ${\rm GL}_2(\mathbb{Q}_p)$, there always will be
infinitely many possibilities for the admissible filtration
(provided there exists at least one). Hence, in some sense, any of
these conjecturally non-zero Banach space representations is
responsible for a specific whole family of Galois representations.

Our picture in the crystalline case is, upon the prize of having and
fixing a square root of $p$ in the coefficient field $K$, very well
adapted to functoriality. The theory of Banach-Hecke algebras as
described above works perfectly well for arbitrary split connected
reductive groups $G$ over $L$ (and is developed in this generality
in section 2), and can in a certain sense be made functorial on the
category of $K$-rational representations of the connected Langlands
dual group $\mathbf{G}'$ of $G$ over $K$. This makes it possible, by
using the formalism of Tannakian categories, to associate with any
specialization of a universal Banach module for $G$ (which still
conjecturally is non-zero) a family of isomorphism classes of
crystalline Galois representations with values in the dual group
$\mathbf{G}'(\overline{K})$ over the algebraic closure
$\overline{K}$ of $K$. Under the already mentioned restrictions, this
was done in \cite{ST} and is here established in general in section
6.

In fact, in \cite{ST} and also in section 6, we have to assume that
the split group $G$ is such that the half sum $\eta$ of its positive
roots is still an integral character. The appearance of this $\eta$
is forced upon us by functoriality. To deal with $\eta$ in general,
we need a counterpart on the Galois side which is a square root of
the cyclotomic character. This does not exist in general on the
Galois group, which leads to the following interesting construction.
The group ${\rm Gal}(\Qpbar/\Qp)$ has, up to isomorphism, a unique
non-trivial central extension ${\rm Gal}(\Qpbar/\Qp)_{(2)}$ with
kernel of order two. Let ${\rm Gal}(\Qpbar/L)_{(2)}$ denote the
restriction of ${\rm Gal}(\Qpbar/\Qp)_{(2)}$ to $L$. On ${\rm
Gal}(\Qpbar/L)_{(2)}$ we always have a square root of the cyclotomic
character provided the coefficient field $K$ is large enough. As a
consequence of the theory of Colmez and Fontaine, one can set up a
theory of crystalline representations of the group ${\rm
Gal}(\Qpbar/L)_{(2)}$. It is exactly families of such which we get
in the case of split groups for which $\eta$ is not integral. This
is the content of section 7. We are grateful to J.-M.\ Fontaine for
helpful discussions about the material in this section.

Throughout the paper we fix two finite extensions $L$ (the base
field) and $K$ (the coefficient field) of $\Qp$ such that
$[L:\Qp]=\vert{\rm Hom}_{\Qp}(L,K)\vert$ where ${\rm
Hom}_{\mathbb{Q}_p}(L,K)$ denotes the set of all
$\mathbb{Q}_p$-linear embeddings of the field $L$ into the field
$K$. We assume $L$ is contained in an algebraic closure $\Qpbar$ of
$\Qp$. We denote by $q=p^f$ the cardinality of the residue field of
$L$ and by $L_0={\rm Frac}(W({\mathbb F}_{q}))$ its maximal
unramified subfield. If $e:=[L:\Qp]/f$, we set $\v(x):=e{\rm
val}_{\Qp}(x)$ (where ${\rm val}_{\Qp}(p):=1$) and
$|x|_L:=q^{-\v(x)}$ for any $x$ in a finite extension of $\Qp$. We
denote by ${\rm W}(\Qpbar/L)$ (resp. $\g$) the Weil (resp. Galois)
group of $L$ and by ${\rm rec}:{\rm W}(\Qpbar/L)^{\rm
ab}\buildrel\sim\over\rightarrow L^{\times}$ the reciprocity map
sending the arithmetic Frobeniuses to the inverse of uniformizers.
If $\lambda\in K^{\times}$, ${\rm unr}(\lambda)$ stands for the
unramified character of $L^{\times}$ sending a uniformizer to
$\lambda$.

\section{Completed Satake-Hecke algebras}\label{satake}

We extend the theory and results of \cite[\S\S 2-3]{ST}.

We fix an $L$-split connected reductive group $\mathbf{G}$ over $L$
and put $G := \mathbf{G}(L)$. Let $(\rho,E)$ be an irreducible
$\mathbb{Q}_p$-rational representation of $G$ in a finite
dimensional $K$-vector space $E$. Fixing a good maximal compact
subgroup $U \subseteq G$, we let $\rho_U := \rho|U$. The
corresponding Satake-Hecke algebra $\mathcal{H}(G,\rho_U)$ is the
convolution algebra over $K$ of all compactly supported functions
$\psi : G \longrightarrow {\rm End}_K(E)$ satisfying:
\begin{equation*}
\psi(u_1gu_2) = \rho(u_1)\circ \psi(g) \circ \rho(u_2)
\end{equation*}
for any $u_1,u_2 \in U$ and $g \in G$. The algebra
$\mathcal{H}(G,\rho_U)$ can naturally be identified with the ring of
$G$-endomorphisms of the compact induction ${\rm ind}^G_U(\rho_U)$.
Fixing once and for all a $U$-invariant norm $\|\ \|$ on $E$, the
$K$-vector space ${\rm ind}^G_U(\rho_U)$ carries the corresponding
$G$-invariant sup-norm also denoted by $\|\ \|$. The $G$-action on
${\rm ind}^G_U(\rho_U)$ extends to an isometric $G$-action on the
completion $B^G_U(\rho_U)$ of ${\rm ind}^G_U(\rho_U)$ with respect
to $\|\ \|$. Using the operator norm on ${\rm End}_K(E)$ we also
have a corresponding sup-norm $\|\ \|$ on $\mathcal{H}(G,\rho_U)$
which is submultiplicative. Its completion $\mathcal{B}(G,\rho_U)$
therefore is a $K$-Banach algebra. It is shown in \cite[Lem.1.3]{ST}
that one has a natural isomorphism of $K$-algebras:
\begin{equation*}
\mathcal{B}(G,\rho_U) \xrightarrow{\ \cong\ }{\rm End}_G^{\rm
cont}(B_U^G(\rho_U))
\end{equation*}
which is an isometry with respect to the operator norm on the right
hand side (which consists of the continuous and $G$-equivariant
endomorphisms of the Banach space $B^G_U(\rho_U)$).

Generalizing the results in \cite[\S\S1-3]{ST}, we want to
explicitly compute the Banach algebra $\mathcal{B}(G,\rho_U)$. In
order to recall what it means for $\rho$ to be
$\mathbb{Q}_p$-rational we introduce the connected reductive group:
\begin{equation*}
\widetilde{\mathbf{G}} := ({\rm Res}_{L/\mathbb{Q}_p} \mathbf{G})_K
\end{equation*}
over $K$ obtained by base extension from the Weil restriction from
$L$ to $\mathbb{Q}_p$ of $\mathbf{G}$. We also put:
\begin{equation*}
\widetilde{G} := \widetilde{\mathbf{G}}(K) = \mathbf{G}(L
\otimes_{\mathbb{Q}_p} K).
\end{equation*}
The ring homomorphism $L \longrightarrow L \otimes_{\mathbb{Q}_p} K$
which sends $a$ to $a \otimes 1$ induces an embedding of groups $G
\hookrightarrow \widetilde{G}$. We have:
\begin{equation*}
\widetilde{\mathbf{G}} = \prod_{\sigma : L\hookrightarrow
K}\mathbf{G}_\sigma
\end{equation*}
where $\mathbf{G}_\sigma$ denotes the base extension of $\mathbf{G}$
to $K$ via the embedding $\sigma : L \rightarrow K$. In particular,
the groups $\mathbf{G}_\sigma$ and $\widetilde{\mathbf{G}}$ are
$K$-split.

The $\mathbb{Q}_p$-rationality of $\rho$ means that there is an
irreducible $K$-rational representation $\widetilde{\rho}$ of
$\widetilde{\mathbf{G}}$ on $E$ such that $\rho$ is the pull-back of
$\widetilde{\rho}$ via $G \hookrightarrow \widetilde{G}$. Since $G$
is Zariski dense in $\widetilde{\mathbf{G}}$ (\cite[\S34.4]{Hu}), the
representation $\widetilde{\rho}$ is uniquely determined by $\rho$.
We also note (\cite[Lem.68]{St}) that:
\begin{equation*}
(\widetilde{\rho},E) \cong \bigotimes_{\sigma : L\hookrightarrow
K}\; (\rho_\sigma,E_\sigma)
\end{equation*}
with irreducible $K$-rational representations
$(\rho_\sigma,E_\sigma)$ of $\mathbf{G}_\sigma$. Conversely any such
tensor product gives rise, by restriction, to an irreducible
$\mathbb{Q}_p$-rational $\rho$.

By a variant of \cite[Lem.1.4]{ST}, one easily shows the following,
where $1_U$ denotes the trivial representation of $U$.

\begin{lem}\label{triv}
The map:
\begin{align*}
 \iota_\rho : \mathcal{H}(G,1_U) & \xrightarrow{\; \cong\;}
   \mathcal{H}(G,\rho_U) \\
 \psi & \longmapsto  \psi\cdot\rho
\end{align*}
is an isomorphism of $K$-algebras.
\end{lem}

At this point we need to introduce further notation. We fix a
maximal $L$-split torus $\mathbf{T}$ in $\mathbf{G}$ and a Borel
subgroup $\mathbf{P} = \mathbf{T}\mathbf{N}$ of $\mathbf{G}$ with
Levi component $\mathbf{T}$ and unipotent radical $\mathbf{N}$. Then
$\widetilde{\mathbf{T}} := ({\rm Res}_{L/\mathbb{Q}_p}
\mathbf{T})_K$ is a maximal $K$-split torus in the Borel subgroup
$\widetilde{\mathbf{P}} := ({\rm Res}_{L/\mathbb{Q}_p}
\mathbf{P})_K$ of $\widetilde{\mathbf{G}}$. We denote by $P$, $T$
and $N$ the group of $L$-valued points of $\mathbf{P}$, $\mathbf{T}$
and $\mathbf{N}$ respectively. The Weyl group of $G$ is the quotient
$W = N(T)/T$ of the normalizer $N(T)$ of $T$ in $G$ by $T$. We
always assume that our fixed maximal compact subgroup $U \subseteq
G$ is special with respect to $T$, and we put $T_0 := U \cap T$ and
$N_0 := U \cap N$. The quotient $\Lambda := T/T_0$ is a free abelian
group of rank equal to the dimension of $\mathbf{T}$ and can
naturally be identified with the cocharacter group
$X_\ast(\mathbf{T})$. Let $\lambda : T \longrightarrow \Lambda$
denote the projection map. The conjugation action of $N(T)$ on $T$
induces $W$-actions on $T$ and $\Lambda$ which we denote by $t
\longmapsto {^w}t$ and $\lambda \longmapsto {^w}\lambda$
respectively. We have the embedding:
\begin{align*}
X^\ast(\mathbf{T}) & \longrightarrow  {\rm Hom}(\Lambda,\mathbb{R}) =: V_{\mathbb{R}} \\
\chi & \longmapsto \v\circ\chi
\end{align*}
which induces an isomorphism:
\begin{equation*}
X^\ast(\mathbf{T}) \otimes \mathbb{R} \xrightarrow{\ \cong\ }
V_{\mathbb{R}}.
\end{equation*}
We therefore may view $V_{\mathbb{R}}$ as the real vector space
underlying the root datum of $G$ with respect to $T$. Evidently any
$\lambda \in \Lambda$ defines a linear form in the dual vector space
$V_{\mathbb{R}}^\ast$ also denoted by $\lambda$. Let $\Phi$ denote
the set of roots of $T$ in $G$ and let $\Phi^+ \subseteq \Phi$ be
the subset of those roots which are positive with respect to $P$. As
usual, $\check{\alpha} \in \Lambda$ denotes the coroot corresponding
to the root $\alpha \in \Phi$. The subset $\Lambda^{--} \subseteq
\Lambda$ of antidominant cocharacters is defined to be the image
$\Lambda^{--} := \lambda(T^{--})$ of:
\begin{equation*}
T^{--} := \{t \in T, |\alpha(t)|_L \geq 1\ \rm{for\ any}\ \alpha \in
\Phi^+\}.
\end{equation*}
Hence, we have:
\begin{equation*}
\Lambda^{--} = \{ \lambda \in \Lambda, \v \circ \alpha(\lambda) \leq
0\ \rm{for\ any}\ \alpha \in \Phi^+ \}.
\end{equation*}

By the Cartan decomposition, $G$ is the disjoint union of the double
cosets $UtU$ with $t$ running over $T^{--}/T_0$. The norm $\|\ \|$
on $\mathcal{H}(G,\rho_U)$ corresponds therefore under the
isomorphism $\iota_\rho$ from Lemma \ref{triv} to the norm $\|\
\|_\rho$ on $\mathcal{H}(G,1_U)$ defined by:
\begin{equation*}
\|\psi\|_\rho := \sup_{t \in T^{--}} (|\psi(t)|_L\cdot\|\rho(t)\|).
\end{equation*}
Hence $\iota_\rho$ extends to an isometric isomorphism of Banach
algebras:
\begin{equation*}
\|\ \|_\rho{\rm -completion\ of}\ \mathcal{H}(G,1_U)
 \xrightarrow{\ \cong\ }
 \mathcal{B}(G,\rho_U).
\end{equation*}
In order to compute this norm further, we let $\widetilde{\xi} \in
X^\ast(\widetilde{\mathbf{T}})$ denote the highest weight (with
respect to $\widetilde{\mathbf{P}}$) of the representation
$\widetilde{\rho}$ and $\xi : T \longrightarrow K^\times$ its
restriction to $T$.

\begin{lem}\label{norm}
We have $\|\rho(t)\| = |\xi(t)|_L$ for any $t \in T^{--}$.
\end{lem}
\begin{proof}
According to \cite[Lem.3.2]{ST}, we have
$\|\widetilde{\rho}(\widetilde{t})\| =
|\widetilde{\xi}(\widetilde{t})|_L$ for any $\widetilde{t} \in
\widetilde{\mathbf{T}}(K)$ which is antidominant with respect to
$\widetilde{\mathbf{P}}$. We therefore have to show that the image
in $\widetilde{\mathbf{T}}(K)$ of any $t \in T^{--}$ is antidominant
with respect to $\widetilde{\mathbf{P}}$. Let $\mathbf{T}_\sigma$,
for any $\sigma \in {\rm Hom}_{\mathbb{Q}_p} (L,K)$, denote the base
extension of $\mathbf{T}$ to $K$ via the embedding $\sigma$. For
simplicity we also write $\sigma : T \hookrightarrow
\mathbf{T}_\sigma(K)$ for the corresponding embedding of groups. The
image of $t$ in $\widetilde{\mathbf{T}}(K) = \prod_\sigma
\mathbf{T}_\sigma(K)$ then is given by $(\sigma(t))_\sigma$. On the
other hand, using the natural identifications $X^\ast(\mathbf{T}) =
X^\ast(\mathbf{T}_\sigma)$, any root $\widetilde{\alpha}$ of
$\widetilde{\mathbf{T}}$ which is positive with respect to
$\widetilde{\mathbf{P}}$ may be viewed as a tuple
$(\alpha_\sigma)_\sigma$ of roots $\alpha_\sigma \in \Phi^+$, and we
have $\widetilde{\alpha}(t) = \prod_\sigma \alpha_\sigma(\sigma(t))
= \prod_\sigma \sigma(\alpha_\sigma(t))$. Since
$|\alpha_\sigma(t)|_L \geq 1$ by assumption, it follows that
$|\widetilde{\alpha}(t))|_L \geq 1$.
\end{proof}

The unnormalized Satake map:
\begin{align*}
 S : \mathcal{H}(G,1_U) & \longrightarrow  K[\Lambda] \\
 \psi & \longmapsto  \sum\limits_{t \in T/T_0}
 \sum\limits_{n \in
 N/N_0} \psi(tn)\lambda(t)
\end{align*}
induces an isomorphism of $K$-algebras:
\begin{equation*}
\mathcal{H}(G,1_U) \xrightarrow{\ \cong\ } K[\Lambda]^{W,\gamma}
\end{equation*}
where the right hand side denotes the $W$-invariants in the group
ring $K[\Lambda]$ with respect to the twisted $W$-action:
\begin{align*}
 W \times K[\Lambda] & \longrightarrow  K[\Lambda] \\
 (w, \sum_\lambda c_\lambda \lambda) & \longmapsto
 w\cdot(\sum_\lambda c_\lambda \lambda) := \sum_\lambda
 \gamma(w,\lambda)c_\lambda {^w\lambda}
\end{align*}
for the $K$-valued cocycle:
\begin{equation*}
\gamma(w,\lambda) := \frac{\delta^{1/2}({^w}\lambda)}
{\delta^{1/2}(\lambda)}
\end{equation*}
with $\delta : P \longrightarrow \mathbb{Q}^\times \subseteq
K^\times$ denoting the modulus character of the Borel subgroup $P$
(compare \cite[\S\S2-3]{ST}). By Lemma \ref{norm} together with a
variant of \cite[Prop.3.5]{ST}, the norm $\|\ \|_\rho$ on
$\mathcal{H}(G,1_U)$ corresponds under this Satake isomorphism to
the restriction of the norm $\|\ \|_\xi$ on $K[\Lambda]$ given by:
\begin{equation*}
\|\sum_{\lambda \in \Lambda} c_\lambda \lambda\|_\xi :=
\sup_{\lambda = \lambda(t)} |\gamma(w,\lambda)\xi({^w}t) c_\lambda|_L
\end{equation*}
with $w \in W$ for each $\lambda$ being chosen in such a way that
${^w}\lambda \in \Lambda^{--}$. One checks that this norm on
$K[\Lambda]$ is submultiplicative and that the twisted $W$-action is
isometric in this norm (compare \cite[Lem.2.1]{ST} and Examples 1
and 2 in \S2 of {\it loc.cit.}). Hence the completion of
$K[\Lambda]$ with respect to $\|\ \|_\xi$ is a $K$-Banach algebra
$K\langle \Lambda;\xi \rangle$ to which the twisted $W$-action
extends. In particular, we may form the Banach algebra $K\langle
\Lambda;\xi \rangle^{W,\gamma}$ of $W$-invariants. As a result of
the discussion so far we obtain the following:

\begin{prop}
The Banach algebras $\mathcal{B}(G,\rho_U)$ and $K\langle
\Lambda;\xi \rangle^{W,\gamma}$ are, in a natural way, isometrically
isomorphic.
\end{prop}

Let $\mathbf{T}'$ denote the $L$-torus dual to $\mathbf{T}$. Its
$K$-valued points are given by $\mathbf{T}'(K) = {\rm
Hom}(\Lambda,K^\times)$. The group ring $K[\Lambda]$ naturally
identifies with the ring of $K$-valued algebraic functions on
$\mathbf{T}'$. We introduce the ``valuation map'':
\begin{equation*}
{\rm val} :\ \mathbf{T}'(K) = {\rm Hom}(\Lambda,K^\times)
\xrightarrow{\ \v\circ\ } {\rm Hom}(\Lambda,\mathbb{R}) =
V_{\mathbb{R}}.
\end{equation*}
In $V_{\mathbb{R}}$ we have the two distinguished points:
\begin{equation*}
\xi_L := \v \circ \xi \qquad{\rm and}\qquad \eta_L :=
[L:\mathbb{Q}_p] \cdot \eta
\end{equation*}
where $\eta$ denotes half the sum of the positive roots in $\Phi^+$.
Let $\leq$ denote the partial order on $V_{\mathbb{R}}$ defined by
$\Phi^+$ (cf. \cite[Chap.VI,\S1.6]{Bo}). Given any point $z \in
V_{\mathbb{R}}$, let $z^{\rm dom}$ be the unique dominant point in
the $W$-orbit of $z$. We put:
\begin{equation*}
V_{\mathbb{R}}^\xi := \{ z \in V_{\mathbb{R}}, (z + \eta_L)^{\rm
dom} \leq \eta_L + \xi_L\}
\end{equation*}
and:
\begin{equation*}
\mathbf{T}'_\xi := {\rm val}^{-1}(V_{\mathbb{R}}^\xi).
\end{equation*}

\begin{thm}\label{affinoid}
\begin{enumerate}
\item[(i)] $V_{\mathbb{R}}^\xi$ is the convex hull of the points
${^w}(\eta_L+\xi_L)- \eta_L$ for $w \in W$;
 \item[(ii)] $\mathbf{T}'_\xi$ is an open $K$-affinoid subdomain of
 the torus $\mathbf{T}'$;
 \item[(iii)] the Banach algebra $K\langle \Lambda;\xi \rangle$ is
naturally isomorphic to the ring of analytic functions on the
affinoid domain $\mathbf{T}'_\xi$;
 \item[(iv)] $K\langle \Lambda;\xi \rangle^{W,\gamma}$ is an affinoid
$K$-algebra.
\end{enumerate}
\end{thm}
\begin{proof}
This is a straightforward variant of \cite[Lem.2.3, Prop.2.4,
Lem.2.7, Ex.3 of \S2]{ST} (compare also the discussion before the
remark in \S6).
\end{proof}

The Weyl group $W$ acts on the affinoid $\mathbf{T}'_\xi$ by:
\begin{equation*}
(w,\zeta) \longrightarrow \frac{{^w}\delta^{1/2}}
{\delta^{1/2}}\cdot{^w}\zeta.
\end{equation*}

\begin{cor}\label{quotientaffinoid}
The Banach algebra $\mathcal{B}(G,\rho_U)$ is naturally isomorphic
to the ring of analytic functions on the quotient affinoid $W
\backslash \mathbf{T}'_\xi$.
\end{cor}

For later purposes (see \S\ref{functoriality}), we also have to
discuss briefly the subsequent renormalization but for which we are
forced to assume that the coefficient field $K$ contains a square
root of $q$. Fix once and for all such a square root $q^{1/2}$. We
then have a preferred square root $\delta^{1/2} \in \mathbf{T}'(K)$
of $\delta \in \mathbf{T}'(K)$ and as a consequence the normalized
Satake isomorphism:
\begin{align*}
S^{\rm norm} : \mathcal{H}(G,1_U) & \xrightarrow{\; \cong\;}
  K[\Lambda]^W \\
 \psi & \longmapsto  \sum\limits_{t \in T/T_0}
 \delta^{-1/2}(t)(\sum\limits_{n \in
 N/N_0} \psi(tn))\lambda(t)
\end{align*}
where now on the right hand side the $W$-invariants are formed with
respect to the action induced by the conjugation action of $N(T)$ on
$T$. We define:
\begin{equation*}
V_{\mathbb{R}}^{\xi,\rm norm} := \{z \in V_{\mathbb{R}}, z^{\rm dom}
\leq \eta_L + \xi_L \}
\end{equation*}
and:
\begin{equation*}
\mathbf{T}'_{\xi,\rm norm} := {\rm val}^{-1}(V_{\mathbb{R}}^{\xi,\rm
norm}).
\end{equation*}
This is an affinoid subdomain which is invariant under the natural
$W$-action on $\mathbf{T}'$. As discussed in \cite{ST} before the
remark in \S6, the above corollary can be reformulated as follows.

\begin{cor}\label{normalized}
The Banach algebra $\mathcal{B}(G,\rho_U)$ is naturally isomorphic
to the ring of analytic functions on the quotient affinoid $W
\backslash \mathbf{T}'_{\xi,\rm norm}$.
\end{cor}

As a consequence of the above corollaries, we have natural
identifications between the set of $K$-valued (continuous)
characters of the Banach algebra $\mathcal{B}(G,\rho_U)$ and the
sets of $K$-rational points $(W \backslash \mathbf{T}'_\xi)(K)$ and
$(W \backslash \mathbf{T}'_{\xi,\rm norm})(K)$ respectively.

\section{Crystalline Galois representations}\label{crystalline}

We focus on the special case $G={\rm GL}_{d+1}(L)$, $d\geq 1$. We
give a link between the constructions and results of \S\ref{satake}
and the theory of crystalline representations, generalizing
\cite[\S5]{ST} (although we do not use the same normalization
as in {\it loc.cit.}).

We let $U:={\rm GL}_{d+1}(\O)$. We fix an irreducible $\Qp$-rational
representation $(\rho,E)$ of $G$ as in \S\ref{satake} and
$(\zeta_1,\cdots,\zeta_{d+1})\in (K^{\times})^{d+1}$. Let $\tau$ be
a permutation of $\{1,\cdots,d+1\}$ such that
$\v(\zeta_{\tau(1)})\leq
\v(\zeta_{\tau(2)})\leq\cdots\leq\v(\zeta_{\tau(d+1)})$. Let
$\widehat\zeta:T\rightarrow K^{\times}$ be the character which sends
$t:={\rm diag}(t_1,\cdots,t_{d+1})\in T$ to
$\prod_j(\zeta_jq^{1-j})^{\v(t_j)}$, in other words:
$$\widehat\zeta:={\rm unr}(\zeta_1)\otimes {\rm unr}(\zeta_2)\abs_L\otimes {\rm unr}(\zeta_3)\abs_L^{2}\otimes\cdots\otimes{\rm unr}(\zeta_{d+1})\abs_L^{d}.$$
With the notations of \S\ref{satake}, we also denote by
$\widehat\zeta$ the $K$-linear map $K[\Lambda]\rightarrow K$ induced
by $\widehat\zeta$ and by $K_{\widehat\zeta}$ the one dimensional
$K$-vector space on which ${\mathcal H}(G,\rho_U)$ acts through the
character:
$${\mathcal H}(G,\rho_U)\buildrel\iota_{\rho}^{-1}\over\longrightarrow {\mathcal H}(G,1_U)\buildrel S\over\longrightarrow K[\Lambda]\buildrel \widehat\zeta\over\longrightarrow K.$$
For $\sigma:L\hookrightarrow K$, we let
$(a_{1,\sigma},\cdots,a_{d+1,\sigma})\in \Z^{d+1}$ with
$a_{j,\sigma}\leq a_{j+1,\sigma}$ be the highest weight of
$\rho_{\sigma}$ with respect to the parabolic subgroup of lower
triangular matrices of $G$ (see \S\ref{satake}).

Recall that an invariant norm on a locally algebraic representation
of $G$ (on a $K$-vector space $E$) is a $p$-adic norm $\|\ \|$ such
that $\|gv\|=\|v\|$ for all $g\in G$ and $v\in E$.

\begin{cor}\label{norm2}
If the locally algebraic representation:
$$K_{\widehat\zeta}\otimes_{{\mathcal H}(G,\rho_U)}{\rm ind}_U^G\rho_U$$
admits an invariant norm, then the following inequalities hold in
$\Q$:
\begin{eqnarray*}
\sum_{j=i}^{d+1}\v(\zeta_{\tau(j)})&\leq &\sum_{j=i}^{d+1}\sum_{\sigma}a_{j,\sigma}+[L:\Qp]{d(d+1)-(i-2)(i-1)\over 2},\ 2\leq i\leq d+1\\
\sum_{j=1}^{d+1}\v(\zeta_{\tau(j)})&=
&\sum_{j=1}^{d+1}\sum_{\sigma}a_{j,\sigma}+[L:\Qp]{d(d+1)\over 2}.
\end{eqnarray*}
\end{cor}
\begin{proof}
If the above locally algebraic representation admits an invariant
norm, then necessarily the image of the unit ball of ${\rm
ind}_U^G\rho_U$ in $K_{\widehat\zeta}\otimes_{{\mathcal
H}(G,\rho_U)}{\rm ind}_U^G\rho_U$ is again a unit ball (i.e. remains
a genuine $\O$-lattice). This implies that the character
$\widehat\zeta:{\mathcal H}(G,\rho_U)\rightarrow K_{\widehat\zeta}$
factors through the completion ${\mathcal B}(G,\rho_U)$ of
${\mathcal H}(G,\rho_U)$, hence defines a $K$-point of ${\bf
T}'_{\xi}$ by Corollary 2.5. By the description of those $K$-points,
this implies that one has, with the notations of \S\ref{satake}:
$$(\v\circ\widehat\zeta + \eta_L)^{\rm dom}\leq \xi_L+\eta_L,$$
that is to say:
\begin{multline*}
\bigg(\big(\v(\zeta_1),\v(\zeta_2)-[L:\Qp],\cdots,\v(\zeta_{d+1})-d[L:\Qp]\big)\\
+[L:\Qp]\Big(-{d\over 2},-{d\over 2}+1,\cdots,{d\over 2}\Big)\bigg)^{\rm dom}\\
\leq
\Big(\sum_{\sigma}a_{1,\sigma},\cdots,\sum_{\sigma}a_{d+1,\sigma}\Big)+[L:\Qp]\Big(-{d\over
2},-{d\over 2}+1,\cdots,{d\over 2}\Big).
\end{multline*}
An immediate computation shows this is equivalent to:
\begin{multline*}
\big(\v(\zeta_1),\v(\zeta_2),\cdots,\v(\zeta_{d+1})\big)^{\rm dom}+[L:\Qp]\Big(-{d\over 2},-{d\over 2},\cdots,-{d\over 2}\Big)\\
\leq
\Big(\sum_{\sigma}a_{1,\sigma},\cdots,\sum_{\sigma}a_{d+1,\sigma}\Big)+[L:\Qp]\Big(-{d\over
2},-{d\over 2}+1,\cdots,{d\over 2}\Big),
\end{multline*}
that is to say:
\begin{multline*}
\big(\v(\zeta_{\tau(1)}),\v(\zeta_{\tau(2)}),\cdots,\v(\zeta_{\tau(d+1)})\big)\leq \Big(\sum_{\sigma}a_{1,\sigma},\cdots,\sum_{\sigma}a_{d+1,\sigma}\Big)\\
+[L:\Qp](0,1,\cdots,d)
\end{multline*}
which is what we want.
\end{proof}

Denote by $\varphi_0:L_0\rightarrow L_0$ the absolute Frobenius. Let
$D:=L_0\otimes_{\Qp}K\cdot e_1 \oplus \cdots \oplus
L_0\otimes_{\Qp}K\cdot e_{d+1}$ be a free $L_0\otimes_{\Qp}K$-module
of rank $d+1$ and denote by $F_{\zeta}$ the unique
$L_0\otimes_{\Qp}K$-linear automorphism of $D$ such that:
$$F_{\zeta}(e_j):=\zeta_j^{-1}e_j,\ \ j\in\{1,\cdots,d+1\}.$$
We call a Frobenius on $D$ any bijective map $\varphi:D\rightarrow
D$ satisfying:
$$\varphi((\ell\otimes k)\cdot d)=(\varphi_0(\ell)\otimes k)\cdot\varphi(d)$$
where $\ell\in L_0$, $k\in K$ and $d\in D$. If $\varphi$ is a
Frobenius on $D$, then $\varphi^f$ is $L_0\otimes_{\Qp}K$-linear.
The isomorphism $L_0\otimes_{\Qp}K\simeq
\prod_{\sigma_0:L_0\hookrightarrow K}K$, $\ell\otimes k\mapsto
(\sigma_0(\ell)k)_{\sigma}$ induces an isomorphism:
$$D\simeq \prod_{\sigma_0:L_0\hookrightarrow K}D_{\sigma_0}$$
where $D_{\sigma_0}:=(0,0,\cdots,0,1_{\sigma_0},0,\cdots,0)\cdot
D_L$. The linear map $\varphi^f$ thus induces a $K$-linear
automorphism on each $D_{\sigma_0}$ and all the pairs
$(D_{\sigma_0},\varphi^f)$ are isomorphic via some power of
$\varphi$. We write $(\varphi^f)^{\rm ss}$ for the semisimple part
of $\varphi^f$ on $D$. We define:
\begin{eqnarray}\label{tn}
t_N(D):={1\over [L:\Qp]}\v\big({\rm
det}_{L_0}(\varphi^f\vert_D)\big).
\end{eqnarray}

Let $D_L:=L\otimes_{L_0}D$, as before the isomorphism
$L\otimes_{\Qp}K\simeq \prod_{\sigma:L\hookrightarrow K}K$,
$\ell\otimes k\mapsto (\sigma(\ell)k)_{\sigma}$ induces an
isomorphism:
$$D_L\simeq \prod_{\sigma:L\hookrightarrow K}D_{L,\sigma}$$
where $D_{L,\sigma}:=(0,0,\cdots,0,1_{\sigma},0,\cdots,0)\cdot D_L$.
To give an $L\otimes_{\Qp}K$-submodule ${\rm Fil}D_L$ of $D_L$ is
thus the same thing as to give a collection $({\rm
Fil}D_{L,\sigma})_{\sigma}$ where ${\rm Fil}D_{L,\sigma}$ is a
$K$-vector subspace of $D_{L,\sigma}$. If $({\rm
Fil}^iD_{L,\sigma})_{i,\sigma}$ is a decreasing exhaustive separated
filtration on $D_L$ by $L\otimes_{\Qp}K$-submodules indexed by $i\in
\Z$, we define:
\begin{eqnarray}\label{th}
t_H(D_L):=\sum_{\sigma}\sum_{i\in \Z}i{\rm dim}_L\big({\rm
Fil}^iD_{L,\sigma}/{\rm Fil}^{i+1}D_{L,\sigma}\big).
\end{eqnarray}
Recall that such a filtration is called admissible (one used to say
weakly admissible) if $t_H(D_L)=t_N(D)$ and if, for any $L_0$-vector
subspace $D'\subseteq D$ preserved by $\varphi$ with the induced
filtration on $D'_L$, one has $t_H(D'_L)\leq t_N(D')$.

For $j\in \{1,\cdots,d+1\}$, let:
\begin{eqnarray}\label{weights}
i_{j,\sigma}:=-a_{d+2-j,\sigma}-(d+1-j)\in \Z.
\end{eqnarray}
Note that one has $i_{1,\sigma}<i_{d,\sigma}<\cdots<i_{d+1,\sigma}$.

\begin{prop}\label{admissible}
The following conditions are equivalent:
\begin{enumerate}
\item[(i)] there is a Frobenius $\varphi$ on $D$ such that $(\varphi^f)^{\rm ss}=F_{\zeta}$ and an admissible filtration $({\rm Fil}^iD_{L,\sigma})_{i,\sigma}$ on the $\varphi$-module $(\varphi,D)$ such that, $\forall\ \sigma$:
$${\rm Fil}^iD_{L,\sigma}/{\rm Fil}^{i+1}D_{L,\sigma}\ne 0\Leftrightarrow i\in\{i_{1,\sigma},\cdots,i_{d+1,\sigma}\};$$
\item[(ii)] the following inequalities hold in $\Q$:
\begin{eqnarray*}
\sum_{j=1}^{i}\sum_{\sigma}i_{j,\sigma}&\leq &-\sum_{j=d+2-i}^{d+1}\v(\zeta_{\tau(j)}),\ 1\leq i\leq d\\
\sum_{j=1}^{d+1}\sum_{\sigma}i_{j,\sigma}&=
&-\sum_{j=1}^{d+1}\v(\zeta_{\tau(j)});
\end{eqnarray*}
\item[(iii)] the (Hodge) polygon associated to: $$\Big(\sum_{\sigma}i_{1,\sigma},\sum_{\sigma}i_{1,\sigma}+\sum_{\sigma}i_{2,\sigma},\cdots,\sum_{j=1}^{d+1}\sum_{\sigma}i_{j,\sigma}\Big)$$
is under the (Newton) polygon associated to:
$$\Big(-\v(\zeta_{\tau(d+1)}),-\v(\zeta_{\tau(d+1)})-\v(\zeta_{\tau(d)}),\cdots,-\sum_{j=1}^{d+1}\v(\zeta_{\tau(j)})\Big)$$
and both have the same endpoints.
\end{enumerate}
\end{prop}
\begin{proof}
(iii) is just a restatement of (ii). Assume (i). By assumption, one
can modify $\tau$ such that the matrix of $\varphi^f$ in the basis
$(e_{\tau(d+2-i)})_{1\leq i\leq d+1}$ is upper triangular with
$(\zeta_{\tau(d+2-i)}^{-1})_i$ on the diagonal. Consider the
subspace $D_{1}:=L_0\otimes_{\Qp}K\cdot e_{\tau(d+1)}$ (which is
preserved by $\varphi$). Viewing $D_1$ (resp. $D_{1,L}$) as just an
$L_0$-vector space (resp. $L$-vector space), one has by (\ref{th})
$\sum_{\sigma}[K:L]i_{1,\sigma}=[K:L]\sum_{\sigma}i_{1,\sigma}\leq
t_H(D_{1,L})$ and by (\ref{tn}) $t_N(D_{1})=[L:\Qp]^{-1}[K:\Qp]{\rm
val}_{L}(\zeta_{\tau(d+1)}^{-1})=-[K:L]\v(\zeta_{\tau(d+1)})$. The
inequality $t_H(D_{1,L})\leq t_N(D_{1})$ then implies the first
inequality $\sum_{\sigma}i_{1,\sigma}\leq -\v(\zeta_{\tau(d+1)})$.
One can proceed with the subspaces:
$$D_i:=L_0\otimes_{\Qp}K\cdot e_{\tau(d+1)}\oplus \cdots\oplus L_0\otimes_{\Qp}K\cdot e_{\tau(d+2-i)}$$
for $2\leq i\leq d+1$ (which are all preserved by $\varphi$). The
inequalities $t_H(D_{i,L})\leq t_N(D_{i})$ for $i\leq d$ imply the
intermediate inequalities of (ii) whereas the equality $t_H(D_{L})=
t_N(D)$ yields the final equality. Assume (ii). Note first that it
is enough to check the admissibility conditions for
$L_0\otimes_{\Qp}K$-submodules preserved by $\varphi$ (instead of
$L_0$-vector subspaces preserved by $\varphi$ but not necessarily
$K$): see \cite[Prop.3.1.1.5]{BM}. Let
$e_{i,\sigma}:=(0,0,\cdots,0,1_{\sigma},0,\cdots,0)\cdot e_i$ for
$\sigma:L\hookrightarrow K$. Define $\varphi$ such that, on each
piece of $D$ where $(\varphi^f)^{\rm ss}$ is scalar, say
$L_0\otimes_{\Qp}K\cdot e_{\tau(i)}\oplus\cdots\oplus
L_0\otimes_{\Qp}K\cdot e_{\tau(i-h)}$, $\varphi^f$ is given by the
matrix in the basis $(e_{\tau(i)},\cdots,e_{\tau(i-h)})$:
$$\begin{pmatrix}
\zeta_{\tau(i)}^{-1}&1&0&\cdots&0\\
0&\ddots&\ddots&\ddots&\vdots\\
\vdots&\ddots&\ddots&\ddots&0\\
\vdots&&\ddots&\ddots&1\\0&\ldots&\ldots&0&\zeta_{\tau(i)}^{-1}
\end{pmatrix}.$$
Define the filtration such that, for each $\sigma$:
\begin{eqnarray*}
{\rm Fil}^iD_{L,\sigma}&:=&D_{L,\sigma}\ {\rm if}\ i\leq i_{1,\sigma}\\
{\rm Fil}^iD_{L,\sigma}&:=&Kf_{j,\sigma}\oplus\cdots\oplus Kf_{d+1,\sigma}\ {\rm if}\ i_{j-1,\sigma}+1\leq i\leq i_{j,\sigma}\\
{\rm Fil}^iD_{L,\sigma}&:=&0\ {\rm if}\ i_{d+1,\sigma}+1\leq i
\end{eqnarray*}
where
$f_{j,\sigma}:=e_{\tau(d+2-j),\sigma}+\lambda_{j,j-1}e_{\tau(d+3-j),\sigma}+\cdots+\lambda_{j,1}e_{\tau(d+1),\sigma}$
if $1\leq j\leq d+1$ and the $\lambda_{j,k}\in K$ are chosen such
that all the determinants for $r<j_1<j_2<\cdots<j_r$:
$$\left| \begin{array}{ccc}
\lambda_{j_1,r} & \cdots & \lambda_{j_1,1}\\
\vdots &\vdots & \vdots\\
\lambda_{j_r,r} & \cdots & \lambda_{j_r,1}
\end{array}\right| $$
are non-zero (which is always generically possible). We then leave
it as an exercice to the reader to check that, on each
$L_0\otimes_{\Qp}K$-submodule $D'$ of $D$ preserved by $\varphi$,
the condition $t_H(D'_L)\leq t_N(D')$ is satisfied (as well as
$t_H(D_L)= t_N(D)$).
\end{proof}

\begin{cor}\label{unsens}
If the locally algebraic representation:
$$K_{\widehat\zeta}\otimes_{{\mathcal H}(G,\rho_U)}{\rm ind}_U^G\rho_U$$
admits an invariant norm, then there is a Frobenius $\varphi$ on $D$
such that $(\varphi^f)^{\rm ss}=F_{\zeta}$ and an admissible
filtration $({\rm Fil}^iD_{L,\sigma})_{i,\sigma}$ on the
$\varphi$-module $(\varphi,D)$ such that, $\forall\ \sigma$, ${\rm
Fil}^iD_{L,\sigma}/{\rm Fil}^{i+1}D_{L,\sigma}\ne 0\Leftrightarrow
i\in\{i_{1,\sigma},\cdots,i_{d+1,\sigma}\}$.
\end{cor}
\begin{proof}
Use (\ref{weights}) to replace the $a_{j,\sigma}$ by the
$i_{j,\sigma}$ in the inequalities of Corollary \ref{norm2} and then
use Proposition \ref{admissible}.
\end{proof}

Using \cite{CF}, we thus get that the existence of an invariant norm
on $K_{\widehat\zeta}\otimes_{{\mathcal H}(G,\rho_U)}{\rm
ind}_U^G\rho_U$ implies the existence of at least one crystalline
representation $V$ of $\g$ of dimension $d+1$ over $K$ such that the
eigenvalues of $\varphi^f$ on $D_{\rm cris}(V):=(B_{\rm
cris}\otimes_{\Qp}V)^{\g}$ are the $\zeta_j^{-1}$ and such that
the Hodge-Tate weights of $V$ are the $-i_{j,\sigma}$.

\section{A general conjecture for de Rham representations}\label{derham}

We keep $G={\rm GL}_{d+1}(L)$, $d\geq 1$. The aim of this section is
to state a conjecture which significantly generalizes and
strengthens Corollary \ref{unsens} (in its statement). We keep the
notations of \S\ref{crystalline}.

Let $L'$ be a finite Galois extension of $L$ and $L'_0$ its maximal
unramified subfield. We assume $[L'_0:\Qp]=\vert{\rm
Hom}_{\Qp}(L'_0,K)\vert$ and we let $p^{f'}$ be the cardinality of
the residue field of $L'_0$ and $\varphi'_0$ be the Frobenius on
$L'_0$ (raising to the $p$ each component of the Witt vectors).
Consider the following two categories:

\begin{enumerate}
\item[(i)]the category ${\rm WD}_{L'/L}$ of representations $(r,N,V)$ of the Weil-Deligne group of $L$ (\cite[\S8]{De}) on a $K$-vector space $V$ of finite dimension such that $r$ is unramified when restricted to ${\rm W}(\Qpbar/L')$;\\

\item[(ii)]the category ${\rm MOD}_{L'/L}$ of quadruples $(\varphi,N,\gL,D)$ where $D$ is a free $L'_0\otimes_{\Qp}K$-modules of finite rank endowed with a Frobenius $\varphi:D\rightarrow D$ as in \S\ref{crystalline}, an $L'_0\otimes_{\Qp}K$-linear endomorphism $N:D\rightarrow D$ such that $N\varphi=p\varphi N$ and an action of $\gL$ commuting with $\varphi$ and $N$ such that $g((\ell\otimes k)\cdot d)=(g(\ell)\otimes k)\cdot g(d)$ ($g\in \gL$, $\ell\in L'_0$, $k\in K$, $d\in D$).
\end{enumerate}

Note that in (ii) $N$ is necessarily nilpotent
(\cite[\S1.1.3]{Fo2}).

There is a functor (due to Fontaine):
$${\rm WD}:{\rm MOD}_{L'/L}\rightarrow {\rm WD}_{L'/L}$$
defined as follows (see \cite{Fo2}). Choose an embedding ${\bf
\sigma}'_0:L'_0\hookrightarrow K$ and let $V:=D_{\sigma'_0}$ (with
$D_{\sigma'_0}$ as in \S\ref{crystalline}). As $N$ is
$L'_0\otimes_{\Qp}K$-linear, it induces a nilpotent $K$-linear
endomorphism again denoted $N:V\rightarrow V$. For $w\in {\rm
W}(\Qpbar/L)$, define $r(w):=\overline w\circ \varphi^{-\alpha(w)}$
where $\overline w$ is the image of $w$ in $\gL$ and $\alpha(w)\in
f\Z$ is the unique integer such that the image of $w$ in ${\rm
Gal}(\Fpbar/\Fp)$ is the $\alpha(w)$-th power of the absolute
arithmetic Frobenius. We immediately see that $r(w)$ is
$L'_0\otimes_{\Qp}K$-linear and thus induces a map again denoted
$r(w):V\rightarrow V$. It is not difficult to check that $(r,N,V)$
is an object of ${\rm WD}_{L'/L}$. Moreover, up to non-natural
isomorphism, the representation $(r,N,V)$ doesn't depend on the
choice of $\sigma'_0$ (see \cite[Lem.2.2.1.2]{BM}).

\begin{prop}\label{equivalence}
The functor ${\rm WD}:{\rm MOD}_{L'/L}\rightarrow {\rm WD}_{L'/L}$
is an equivalence of categories.
\end{prop}
\begin{proof}
We build a quasi-inverse. Let $(r,N,V)$ be an object of ${\rm
WD}_{L'/L}$. As an $L'_0\otimes_{\Qp}K$-module, we take:
$$D:=\bigoplus_{n=0}^{f'-1}V_{\sigma'_0\circ{\varphi'_0}^{-n}}$$
where $V_{\sigma'_0\circ{\varphi'_0}^{-n}}=V$ but with $L'_0$ acting
via $\sigma'_0\circ{\varphi'_0}^{-n}$. We define
$\varphi:D\rightarrow D$ by mapping identically
$V_{\sigma'_0\circ{\varphi'_0}^{-n}}$ to
$V_{\sigma'_0\circ{\varphi'_0}^{-n-1}}$ if $0\leq n<f'-2$ and by
mapping $V_{\sigma'_0\circ{\varphi'_0}^{1-f}}$ to $V_{\sigma'_0}$ by
any geometric Frobenius of ${\rm W}(\Qpbar/L')$. We define
$N:V_{\sigma'_0}\rightarrow V_{\sigma'_0}$ as the endomorphism $N$
on $V$ and $N:V_{\sigma'_0\circ{\varphi'_0}^{-n}}\rightarrow
V_{\sigma'_0\circ{\varphi'_0}^{-n}}$ for $1\leq n\leq f'-1$ as
$p^n\varphi^n\circ N\circ\varphi^{-n}$. One checks that
$N\varphi=p\varphi N$ on $D$. For any $g\in \gL$, let $w\in {\rm
W}(\Qpbar/L)$ be a lifting of $g$ and define
$g:V_{\sigma'_0\circ{\varphi'_0}^{-n}}\rightarrow
V_{\sigma'_0\circ{\varphi'_0}^{-n-\alpha(w)}}$ for $0\leq n\leq
f'-1$ as $r(w)\circ \varphi^{\alpha(w)}$ (where
$r(w):V_{\sigma'_0\circ{\varphi'_0}^{-n-\alpha(w)}}\rightarrow
V_{\sigma'_0\circ{\varphi'_0}^{-n-\alpha(w)}}$ is the action of
$w\in {\rm W}(\Qpbar/L)$ on $V$). As $r\mid_{{\rm W}(\Qpbar/L')}$ is
unramified, one checks that this doesn't depend on the choice of the
lifting $w$ and that the action of $g$ commutes with $\varphi$ and
$N$. The final details are left to the reader.
\end{proof}

If $D$ is an object of ${\rm MOD}_{L'/L}$, we define similarly to
(\ref{tn}):
\begin{eqnarray}\label{tng}
t_N(D):={1\over[L:L_0]f'}\v\big({\rm
det}_{L'_0}(\varphi^{f'}\vert_D)\big).
\end{eqnarray}
For $\sigma:L\hookrightarrow K$, let:
$$D_{L',\sigma}:=D_{L'}\otimes_{L'\otimes_{\Qp}K}(L'\otimes_{L,\sigma}K).$$
One has again $D_{L'}\simeq \prod_{\sigma:L\hookrightarrow
K}D_{L',\sigma}$. To give an $L'\otimes_{\Qp}K$-submodule ${\rm
Fil}D_{L'}$ of $D_{L'}$ preserved by $\gL$ is thus the same thing as
to give a collection $({\rm Fil}D_{L',\sigma})_{\sigma}$ where ${\rm
Fil}D_{L',\sigma}$ is a free $L'\otimes_{L,\sigma}K$-submodule of
$D_{L',\sigma}$ (hence a direct factor as
$L'\otimes_{L,\sigma}K$-modules) preserved by the action of $\gL$.
If $({\rm Fil}^iD_{L',\sigma})_{i,\sigma}$ is a decreasing
exhaustive separated filtration on $D_{L'}$ by
$L'\otimes_{\Qp}K$-submodules indexed by $i\in \Z$ and preserved by
$\gL$, we define similarly to (\ref{th}):
\begin{eqnarray}\label{thg}
t_H(D_{L'}):=\sum_{\sigma}\sum_{i\in \Z}i{\rm dim}_{L'}\big({\rm
Fil}^iD_{L',\sigma}/{\rm Fil}^{i+1}D_{L',\sigma}\big).
\end{eqnarray}
Recall that such a filtration is called admissible if
$t_H(D_L)=t_N(D)$ and if, for any $L'_0$-vector subspace
$D'\subseteq D$ preserved by $\varphi$ and $N$ with the induced
filtration on $D'_{L'}$, one has $t_H(D'_{L'})\leq t_N(D')$.

Fix a choice of $q^{1/2}$ in $\Qpbar$. If $(r,N,V)$ is an object of
${\rm WD}_{L'/L}$ such that $r$ is semisimple, we denote by
$\pi^{\rm unit}$ the smooth irreducible representation of $G$ over
$\Qpbar$ corresponding to $(r,N,V)$ by the unitary local Langlands
correspondence normalized so that the central character of $\pi^{\rm
unit}$ is ${\rm det}(r,N,V)\circ{\rm rec}^{-1}$. Note that $\pi^{\rm
unit}$ depends on the choice of $q^{1/2}$.

We now modify the unitary local Langlands correspondence as follows.

Assume first that $\pi^{\rm unit}$ is generic (\cite[\S2.3]{Ku}).
The representation:
$$\pi^{\rm unit}\otimes_{\Qpbar}|{\rm det}|_L^{-d/2}$$
is the extension of scalars from $K$ to $\Qpbar$ of a unique
irreducible admissible representation over $K$ which doesn't depend
any-more on the choice of $q^{1/2}$ (use \cite[Conj.4.4]{Cl} proved
in \cite[\S7]{He} and \cite[Prop.3.2]{Cl}). Denote by $\pi$ this
irreducible generic representation over $K$.

Assume now that $\pi^{\rm unit}$ is not generic. The Langlands
classification tells us that $\pi^{\rm unit}$ is the unique quotient
of a normalized parabolic induction:
\begin{equation}\label{parabolic}
{\rm Ind}_{Q}^GL(b_1,\tau_1)\otimes\cdots\otimes L(b_s,\tau_s)
\end{equation}
where the $\tau_i$ are irreducible supercuspidal representations of
${\rm GL}_{n_i}(L)$, the $b_i$ are positive integers, the
$L(b_i,\tau_i)$ are the associated generalized Steinberg (same
notation as in \cite[\S3.1]{Cl}) and $Q$ is the upper parabolic
subgroup of $G$ of Levi subgroup isomorphic to ${\rm
GL}_{b_1n_1}(L)\times\cdots\times {\rm GL}_{b_sn_s}(L)$. In fact,
writing $(r,N,V)=\oplus_i(r_i,N_i,V_i)$ over $\Qpbar$ where
$(r_i,N_i,V_i)$ is indecomposable, $L(b_i,\tau_i)$ corresponds to
$(r_i,N_i,V_i)$ by the above unitary local Langlands correspondence.
In (\ref{parabolic}), the $L(b_i,\tau_i)$ are ordered so that the
``does not precede'' condition of \cite[Def.1.2.4]{Ku} holds and the
representation (\ref{parabolic}) doesn't depend on such an order
(the proof of this fact is the same as that of \cite[Prop.6.4]{Ze}
using \cite[Th.9.7(a)]{Ze} instead of \cite[Th.4.2]{Ze}).

\begin{lem}\label{model}
The representation:
$$\big({\rm Ind}_{Q}^GL(b_1,\tau_1)\otimes\cdots\otimes L(b_s,\tau_s)\big)\otimes_{\Qpbar}|{\rm det}|_L^{-d/2}$$
admits a unique model over $K$ which doesn't depend on the choice of
$q^{1/2}$.
\end{lem}
\begin{proof}
Let ${\mathcal L}(b_i,\tau_i):=L(b_i,\tau_i)\otimes_{\Qpbar}|{\rm
det}|_L^{(1-b_in_i)/2}$, then ${\mathcal L}(b_i,\tau_i)$ doesn't
depend on the choice of $q^{1/2}$ and one has:
\begin{multline}\label{parabolicnorm}
{\rm Ind}_{Q}^GL(b_1,\tau_1)\otimes\cdots\otimes L(b_s,\tau_s)\otimes_{\Qpbar}|{\rm det}|_L^{-d/2}=\\
{\rm ind}_{Q}^G{\mathcal L}(b_1,\tau_1)\otimes{\mathcal
L}(b_2,\tau_2)|{\rm det}|_L^{-b_1n_1}\otimes\cdots\otimes {\mathcal
L}(b_s,\tau_s)|{\rm det}|_L^{-\sum_{j=1}^{s-1}b_jn_j}
\end{multline}
where the induction on the right hand side is an unnormalized
parabolic induction (no twist by any modulus). The group ${\rm
Gal}(\Qpbar/K)$ acts on the set of isomorphism classes of smooth
representations of $G$ over $\Qpbar$ by sending a representation to
the class of its twist by an element of ${\rm Gal}(\Qpbar/K)$ (see
\cite[\S3.1]{Cl}). As $r$ is a $K$-representation and as the
correspondence $(r,N,V)\mapsto \pi^{\rm
unit}\otimes_{\overline\Qp}|{\rm det}|_L^{-d/2}$ commutes with the
action of ${\rm Gal}(\Qpbar/K)$ (\cite[\S7]{He}), one has that ${\rm
Gal}(\Qpbar/K)$ permutes the representations ${\mathcal
L}(b_i,\tau_i)$ without changing the values of $b_i$ and $n_i$. One
can then reorder the ${\mathcal L}(b_i,\tau_i)$ in the right hand
side of (\ref{parabolicnorm}) according to the orbits of this action
of ${\rm Gal}(\Qpbar/K)$ and rewrite the parabolic induction as a
step by step parabolic induction, each inducing representation
corresponding to an orbit. Using \cite[Th.9.7]{Ze}, one checks that
each inducing representation is irreducible and fixed by ${\rm
Gal}(\Qpbar/K)$. By \cite[Prop.3.2]{Cl}, each such inducing
representation is then the extension of scalars to $\Qpbar$ of a
unique model defined over $K$. Inducing these models, one gets like
this a model over $K$ of the representation (\ref{parabolicnorm}).
Finally, using the Schur lemma for the representation
(\ref{parabolicnorm}) (which holds because it is of finite length
and has a unique irreducible quotient which occurs with multiplicity
$1$, see \cite[\S\S7-9]{Ze}), the same proof as in
\cite[Prop.3.2]{Cl} shows that this model over $K$ is unique. It
doesn't depend on the choice of $q^{1/2}$ as no representation in
this proof does.
\end{proof}

We call $\pi$ the unique model over $K$ given by Lemma \ref{model}.

If $(r,N,V)$ is an object of ${\rm WD}_{L'/L}$, we denote by
$(r,N,V)^{\rm ss}\in {\rm WD}_{L'/L}$ its $F$-semisimplification
(i.e. the underlying Weil representation is the
semisimplifica\-tion of $r$, see \cite[\S8.5]{De}).

We now fix:

\begin{enumerate}
\item[(i)]an object $(r,N,V)$ of ${\rm WD}_{L'/L}$ such that $r$ is semisimple;\\

\item[(ii)]for each $\sigma:L\hookrightarrow K$, a list of $d+1$ integers $i_{1,\sigma}<\cdots<i_{d+1,\sigma}$.
\end{enumerate}

From (i), we define as above the smooth admissible representation
$\pi$. From (ii), we define for $\sigma:L\hookrightarrow K$ and
$j\in \{1,\cdots,d+1\}$:
\begin{eqnarray}\label{weightsinverse}
a_{j,\sigma}:=-i_{d+2-j,\sigma}-(j-1)
\end{eqnarray}
(note that $a_{1,\sigma}\leq a_{2,\sigma}\leq\cdots\leq
a_{d+1,\sigma}$) and we denote by $\rho$ the unique $\Qp$-rational
representation of $G$ over $K$ such that
$\widetilde\rho=\otimes_{\sigma}\rho_{\sigma}$ with $\rho_{\sigma}$
of highest weight $(a_{1,\sigma},\cdots,a_{d+1,\sigma})$ (see
\S\ref{satake}).

We now state our main conjecture:

\begin{conj}\label{conjecture}
Fix $(r,N,V)$ as in (i), $(i_{j,\sigma})_{j,\sigma}$ as in (ii) and
define $\pi$ and $\rho$ as above. The following conditions are
equivalent:
\begin{enumerate}
\item[(i)] there is an invariant norm on $\rho\otimes_K\pi$;\\
\item[(ii)] there is an object $(\varphi,N,\gL,D)$ of ${\rm MOD}_{L'/L}$ such that:
$${\rm WD}\big(\varphi,N,\gL,D\big)^{\rm ss}=(r,N,V)$$
and an admissible filtration $({\rm Fil}^iD_{L',\sigma})_{i,\sigma}$
preserved by $\gL$ on $D_{L'}$ such that, $\forall\ \sigma$:
$${\rm Fil}^iD_{L',\sigma}/{\rm Fil}^{i+1}D_{L',\sigma}\ne 0\Leftrightarrow i\in\{i_{1,\sigma},\cdots,i_{d+1,\sigma}\}.$$
\end{enumerate}
\end{conj}

Using \cite{CF}, this conjecture predicts that the existence of an
invariant norm on $\rho\otimes_K\pi$ is {\it equivalent} to the
existence of a potentially semi-stable representation $V$ of $\g$
that has dimension $d+1$ over $K$, such that its Hodge-Tate weights
are the $-i_{j,\sigma}$ and such that the $F$-semisimplification of
its associated Weil-Deligne representation (\cite[\S2.3.7]{Fo2}) has
$\pi$ as Langlands parameter (modified as above).

\begin{rem}
{\rm Note that we {\it do not} require the stronger statement that
${\rm WD}(\varphi,N,\gL,D)\!=(r,N,V)$ in Conjecture
\ref{conjecture}. Indeed, the conjecture would be false in that
case: take $d=1$, $L=L'=\Qp$, $N=0$ and $r$ scalar, then there is no
module $(\varphi,D)$ with an admissible filtration such that ${\rm
WD}(\varphi,D)=(r,V)$ (as $\varphi$ is scalar), but, at least for
small weights, there is an invariant norm on $\rho\otimes_K\pi$
(\cite[Th.1.3]{Br}).}
\end{rem}

\begin{rem}
{\rm Replacing $L'_0$ by the maximal unramified extension $\Qp^{\rm
nr}$ of $\Qp$ and assuming that $K$ is a finite extension of the
$p$-adic completion of $\Qp^{\rm nr}$, one can state an equivalence
of categories analogous to that of Proposition \ref{equivalence}
without having to specify $L'$. However, it doesn't seem to be
written in the literature that any admissible filtration preserved
by $\gL$ on $D_{L'}$ still corresponds to a $p$-adic semi-stable
representation of $\g$ with the right dimension over $K$. Although
it should be probably easy to derive such a proof from \cite{CF}, we
have preferred to limit ourselves to the statement as in
\ref{conjecture}, depending on some $L'$.}
\end{rem}

Assume for simplicity that the Jordan-H\"older constituents (over
$\Qpbar$) of the Weil representation $r$ are pairwise
non-isomorphic, so that there is only one object $(\varphi,N,\gL,D)$
of ${\rm MOD}_{L'/L}$ such that ${\rm
WD}\big(\varphi,N,\gL,D\big)^{\rm ss}=(r,N,V)$, namely the one given
by Proposition \ref{equivalence}. Then a natural question suggested
by Conjecture \ref{conjecture} would be to ask if a {\it given}
invariant norm on $\rho\otimes_K\pi$ corresponds to a {\it specific}
admissible filtration on $D_{L'}$, and conversely. This seems
(roughly) to hold at least in the case $G={\rm GL}_2(\Qp)$, but we
lack examples so far for more general cases.

\section{Partial cases of the conjecture}

We keep the notations of \S\ref{derham}. We check here several
special cases of (weak forms of) Conjecture \ref{conjecture}. We
start with an observation on the central character of
$\rho\otimes_K\pi$.

\begin{prop}\label{central}
We keep the same notations as in Conjecture \ref{conjecture}. The
following conditions are equivalent:
\begin{enumerate}
\item[(i)] there is an invariant norm on the central character of $\rho\otimes_K\pi$;\\
\item[(ii)] for any (or equivalently one) object $(\varphi,N,\gL,D)$ of ${\rm MOD}_{L'/L}$ such that ${\rm WD}(\varphi,N,\gL,D)^{\rm ss}=(r,N,V)$ and any (or equivalently one) filtration $({\rm Fil}^iD_{L',\sigma})_{i,\sigma}$ preserved by $\gL$ on $D_{L'}$ such that, $\forall\ \sigma$:
$${\rm Fil}^iD_{L',\sigma}/{\rm Fil}^{i+1}D_{L',\sigma}\ne 0\Leftrightarrow i\in\{i_{1,\sigma},\cdots,i_{d+1,\sigma}\},$$
one has $t_H(D_{L'})=t_N(D)$.
\end{enumerate}
\end{prop}
\begin{proof}
Let $\chi_{\rho}$ (resp. $\chi_{\pi}$) be the central character of
$\rho$ (resp. $\pi$). There is an invariant norm on the central
character of $\rho\otimes_K\pi$ if and only if:
\begin{eqnarray}\label{integral}
\v(\chi_{\rho})(\pi_L)+\v(\chi_{\pi})(\pi_L)=0.
\end{eqnarray}
Choose an embedding $\sigma'_0:L'_0\hookrightarrow K$, one has:
\begin{eqnarray*}
\v(\chi_{\rho}(\pi_L)) &=& \sum_{\sigma}\sum_{j=1}^{d+1}a_{j,\sigma}\\
&=&\sum_{j=1}^{d+1}\sum_{\sigma}a_{j,\sigma}\\
\v(\chi_{\pi}(\pi_L)) &=& -\v\big(({\rm det}_K(r))({\rm Frob.\ arith.})\big)+[L:\Qp]{d(d+1)\over 2}\\
&=&{f\over f'}\v\big({\rm
det}_K(\varphi^{f'}\vert_{D_{\sigma'_0}})\big)+[L:\Qp]{d(d+1)\over
2}.
\end{eqnarray*}
On the other hand, for any object $(\varphi,N,\gL,D)$ and any
filtra\-tion as in the statement, one has by (\ref{thg}) and
(\ref{tng}):
\begin{eqnarray*}
t_H(D_{L'}) &=&\sum_{\sigma}\sum_{j=1}^{d+1}[K:L]i_{j,\sigma}\\
&\buildrel{(\ref{weightsinverse})}\over =&-[K:L]\sum_{j=1}^{d+1}\sum_{\sigma}a_{j,\sigma}-[K:\Qp]{d(d+1)\over 2}\\
t_N(D) &=& {1\over[L:L_0]f'}\v\big({\rm det}_{L'_0}(\varphi^{f'}\vert_D)\big)\\
&=&{[K:\Qp]\over [L:L_0]f'}\v\big({\rm det}_K(\varphi^{f'}\vert_{D_{\sigma'_0}})\big)\\
&=&[K:L]{f\over f'}\v\big({\rm
det}_K(\varphi^{f'}\vert_{D_{\sigma'_0}})\big).
\end{eqnarray*}
One sees that:
$$t_N(D)-t_H(D_{L'})=[K:L]\big(\v(\chi_{\rho})(\pi_L)+\v(\chi_{\pi})(\pi_L)\big)$$
which implies the proposition by (\ref{integral}).
\end{proof}

We first look at the supercuspidal case of Conjecture
\ref{conjecture} which turns out to be easy.

\begin{thm}\label{supercuspidal}
When $r$ is absolutely irreducible, Conjecture \ref{conjecture} is
true.
\end{thm}
\begin{proof}
In that case, $\pi$ is a supercuspidal representation. Hence $\pi$,
and thus also $\rho\otimes\pi$, can be written as compact inductions
from a compact open subgroup modulo the center. It is easily checked
that such a compact induction admits an invariant norm if and only
if its central character admits an invariant norm. On the Galois
side, it is enough to check the admissibility conditions for
$L'_0\otimes_{\Qp}K$-submodules of $D$ preserved by $\varphi$ and
$\gL$ (see \cite[Prop.3.1.1.5]{BM} and \cite[Prop.4.4.9]{Fo1}). But
because of Proposition \ref{equivalence} and the assumption, the
only non-zero such module is $D$ itself, hence the admissibility
conditions fall down to just $t_H(D_{L'})=t_N(D)$. The result then
follows from Proposition \ref{central}, using
$D_{L',\sigma}=L'\otimes_L(D_{L',\sigma})^{\gL}$ (Hilbert 90) to
build free $L'\otimes_{L,\sigma}K$-submodules ${\rm
Fil}^iD_{L',\sigma}\subseteq D_{L',\sigma}$ preserved by $\gL$.
\end{proof}

We then turn to the more general case of generalized Steinberg
representations.

\begin{prop}\label{steinberg}
With $(r,N,V)$, $(i_j)_{j,\sigma}$, $\rho\otimes_K\pi$ as in
Conjecture \ref{conjecture} assume that $(r,N,V)$ is indecomposable
over $\Qpbar$ but not irreducible. Let $(r,N,V) = {\rm
WD}\big(\varphi,N,\gL,D\big)$ with $(\varphi,N,\gL,D)$ in ${\rm
MOD}_{L'/L}$. Then the following conditions are equivalent:
\begin{enumerate}
\item[(i)] there is an invariant norm on the central character of $\rho\otimes_K\pi$;\\
\item[(ii)] there is an admissible filtration $({\rm Fil}^iD_{L',\sigma})_{i,\sigma}$ preserved by $\gL$ on $D_{L'}$ such that, $\forall\ \sigma$:
$${\rm Fil}^iD_{L',\sigma}/{\rm Fil}^{i+1}D_{L',\sigma}\ne 0\Leftrightarrow i\in\{i_{1,\sigma},\cdots,i_{d+1,\sigma}\}.$$
\end{enumerate}
\end{prop}
\begin{proof}
Note first that the condition ${\rm
WD}\big(\varphi,N,\gL,D\big)=(r,N,V)$ is here equivalent to the
condition ${\rm WD}\big(\varphi,N,\gL,D\big)^{\rm ss}=(r,N,V)$. We
have already seen in Proposition \ref{central} that (ii) implies
(i). Assume (i) and let $D_0:={\rm Ker}(N:D\rightarrow D)$ which is,
by assumption and using Proposition \ref{equivalence}, a simple
object of ${\rm MOD}_{L'/L}$ (with the induced actions of $\varphi$
and $\gL$). The assumptions on $r$ imply that we necessarily have
$D_0\subsetneq D$. Let $d_0+1$ be the rank of $D_0$ over
$L'_0\otimes_{\Qp}K$ and $s\in {\mathbb N}$ such that
$(s+1)(d_0+1)=d+1$. In ${\rm MOD}_{L'/L}$, $(\varphi,N,\gL,D)$ can
be described as $D_0\oplus D_0(1)\oplus \cdots\oplus D_0(s)$ where
$D_0(n):=D_0$ but with $\varphi$ multiplied by $p^n$ (and same
action of $\gL$) and where $N:D\rightarrow D$ is $0$ on $D_0$ and
sends $D_0(n)$ to $D_0(n-1)$ by the identity map if $n>0$. The only
subobjects of $D$ in ${\rm MOD}_{L'/L}$ are thus $D_0\oplus
D_0(1)\oplus \cdots\oplus D_0(n)$ for $0\leq n\leq s$. Let $({\rm
Fil}^iD_{L',\sigma})_{i,\sigma}$ be any decreasing separated
exhaustive filtration on $D_{L'}$ preserved by $\gL$ such that, for
$0\leq j\leq s$ and all $\sigma:L\hookrightarrow K$:
$${\rm Fil}^{i_{j(d_0+1)+1,\sigma}}D_{L',\sigma}:=D_0(j)_{L',\sigma}\oplus\cdots\oplus D_0(s)_{L',\sigma}.$$
Then it is not difficult to check that, for $0\leq i\leq s$:
\begin{multline*}
t_N(D_0\oplus D_0(1)\oplus \cdots\oplus D_0(n))=[K:L]\Big((n+1){t_N(D_0)\over [K:L]}\\
+[L:\Qp](d_0+1)(1+2+\cdots+n)\Big)
\end{multline*}
and that by (\ref{thg}):
\begin{multline*}
t_H((D_0\oplus D_0(1)\oplus \cdots\oplus D_0(n))_{L'})=[K:L]\bigg(\sum_{j=1}^{d_0+1}\sum_{\sigma}i_{j,\sigma}+\sum_{j=d_0+2}^{2d_0+2}\sum_{\sigma}i_{j,\sigma}+\cdots\\
\cdots+\sum_{j=n(d_0+1)+1}^{(n+1)(d_0+1)}\sum_{\sigma}i_{j,\sigma}\bigg).
\end{multline*}
Applying Lemma \ref{combi} below with:
\begin{eqnarray*}
i_n&:=&\sum_{j=n(d_0+1)+1}^{(n+1)(d_0+1)}\sum_{\sigma}i_{j,\sigma}\\
c&:=&[K:L]^{-1}t_N(D_0)\\
h&:=&[L:\Qp](d_0+1)
\end{eqnarray*}
yields the inequalities for $0\leq n\leq s$:
$$t_H((D_0\oplus D_0(1)\oplus \cdots\oplus D_0(n))_{L'})\leq t_N(D_0\oplus D_0(1)\oplus \cdots\oplus D_0(n))$$
(the last one being an equality) which exactly mean that the
filtration $({\rm Fil}^iD_{L',\sigma})_{i,\sigma}$ is admissible.
\end{proof}

\begin{lem}\label{combi}
Let $s\in {\mathbb N}$, $i_0,\cdots,i_s\in \Z$, $h\in \Z$ and $c\in
\Q$. Assume that:
\begin{eqnarray}
\label{ecart}i_{n-1}+h&\leq &i_n\ {\rm for}\ 1\leq n\leq s\\
\label{sumtos}i_0+\cdots+i_s &\leq & (s+1)c+h(1+2+\cdots+s).
\end{eqnarray}
Then, for $0\leq n\leq s$, one has the inequalities:
$$i_0+\cdots+i_n\leq (n+1)c+h(1+2+\cdots +n).$$
\end{lem}
\begin{proof}
From (\ref{ecart}), we get:
$$i_0+\cdots+i_s\geq i_0+(i_0+h)+\cdots +(i_0+sh)=(s+1)i_0+h(1+2+\cdots +s).$$
Using (\ref{sumtos}), we immediately deduce $c\geq i_0$ (case $n=0$
of the above inequalities). Assume that, for some $n\in
\{1,\cdots,s\}$, we have:
\begin{eqnarray}\label{nmin}
i_0+\cdots+i_n> (n+1)c+h(1+2+\cdots +n)
\end{eqnarray}
and choose the smallest such $n$. As $nc+h(1+2+\cdots +n-1)\geq
i_0+\cdots+i_{n-1}$, we have:
$$i_0+\cdots+i_n> i_0+\cdots+i_{n-1}+c+hn,$$
hence $i_n> c+hn$. Using (\ref{ecart}), we thus get $i_{n+1}>
c+h(n+1),\cdots,i_{s}> c+hs$, hence:
\begin{eqnarray}\label{nmax}
i_{n+1}+\cdots +i_s&> &c+h(n+1)+\cdots +c+hs\\
\nonumber&=&(s-n)c+h(n+1+\cdots +s).
\end{eqnarray}
Adding (\ref{nmin}) and (\ref{nmax}), we get:
$$i_0+\cdots +i_s> (s+1)c+h(1+2+\cdots s)$$
which is in contradiction with (\ref{sumtos}).
\end{proof}

Combining Proposition \ref{steinberg} with Conjecture
\ref{conjecture}, we get in particular the following conjecture:

\begin{conj}\label{generalized}
Let $\rho$ be an irreducible algebraic $\Qp$-rational representation
of $G$ over $K$ and let $\pi$ be a generalized Steinberg
representation of $G$ over $K$. Then $\rho\otimes_K\pi$ admits an
invariant norm if and only if its central character is integral.
\end{conj}

To deduce Conjecture \ref{generalized} from Conjecture
\ref{conjecture}, note that $(r,N,V)$ is indecomposable over
$\Qpbar$ if and only if $\pi$ is a generalized Steinberg (and thus
always generic), that $\rho\otimes_K\pi$ admits an invariant norm if
and only if $\rho\otimes_K\pi\otimes_KK'$ does (where $K'$ is a
finite extension of $\Qp$) and then use Proposition \ref{steinberg}
(replacing $K$ by some $K'$ containing $L'_0$ if necessary).

Finally, we look at the case of unramified principal series.

\begin{thm}\label{conjunramified}
Assume $(r,N,V)$ is an unramified $K$-split Weil representation,
then (i) implies (ii) in Conjecture \ref{conjecture}.
\end{thm}
\begin{proof}
By assumption $N=0$, $L'=L$ and $(r,V)$ sends any arithmetic
Frobenius of ${\rm W}(\Qpbar/L)$ to ${\rm
diag}(\zeta_1,\cdots,\zeta_{d+1})$ for some
$(\zeta_1,\cdots,\zeta_{d+1})\in ({K}^{\times})^{d+1}$. Let
$\widehat\zeta$ be as in \S\ref{crystalline}. Then, an examination
of the proof of  \cite[Lem.3.1]{Da} shows that one has
$K_{\widehat\zeta}\otimes_{{\mathcal H}(G,\rho_U)}{\rm
ind}_U^G\rho_U$ isomorphic to $\rho\otimes_K\pi$. Hence (i) is
equivalent to the existence of an invariant norm on
$K_{\widehat\zeta}\otimes_{{\mathcal H}(G,\rho_U)}{\rm
ind}_U^G\rho_U$. Going back through the definition of the functor
${\rm WD}$, one sees that the result then exactly follows from
Corollary \ref{unsens}.
\end{proof}

When $(r,N,V)$ is a Weil representation that is a direct sum of
characters, the (ii) $\Rightarrow$ (i) sense in Conjecture
\ref{conjecture} seems much deeper (even if the characters are
unramified). The only known case so far is $d=1$, $L=\Qp$ and $r$
non-scalar up to twist (\cite{BB}).

\begin{rem}
{\rm As in Proposition \ref{admissible}, starting from $(r,N,V)\in
{\rm WD}_{L'/L}$, one can give explicit conditions which are
equivalent to the existence of an object $(\varphi,N,\gL,D)$
equipped with an admissible filtration as in Conjecture
\ref{conjecture}, at least for $K$ large enough. Assume one can
write:
$$(r,N,V)=\oplus_{i\in \{1,\cdots,s\}}(r_i,N_i,V_i)$$
over $K$ with $(r_i,N_i,V_i)$ absolutely indecomposable and denote
simply by $D_i$ the object of ${\rm MOD}_{L'/L}$ such that ${\rm
WD}(D_i)=(r_i,N_i,V_i)$. Let $t_{N,i}:=t_N(D_i)$ as in (\ref{tng})
and $d_i:={\rm dim}_K(r_i)$. Order the set of representations
$(r_i,N_i,V_i)$ so that $t_{N,i}$ increases when $i$ grows, and
inside each subset where $t_{N,i}$ is constant so that $d_i$
decreases with $i$ grows. Then, as in Proposition \ref{admissible},
one can prove using Lemma \ref{combi} that there exists
an object $(\varphi,N,\gL,D)$ of ${\rm MOD}_{L'/L}$ equipped with an admissible filtration as in Conjecture \ref{conjecture} if and only if the polygon associated
to:
$$\Big(\sum_{j=1}^{d_1}\sum_{\sigma}i_{j,\sigma},\sum_{j=1}^{d_1+d_2}\sum_{\sigma}i_{j,\sigma},\cdots,\sum_{j=1}^{d+1}\sum_{\sigma}i_{j,\sigma}\Big)$$
is under the polygon associated to
$\big(t_{N,1},t_{N,1}+t_{N,2},\cdots,\sum_{j=1}^st_{N,j}\big)$ and
both have the same endpoints. Using Conjecture \ref{conjecture},
this gives a conjectural necessary and sufficient explicit condition
for $\rho\otimes_K\pi$ to admit an invariant norm.}
\end{rem}

\section{Towards a $p$-adic unramified functoriality I}\label{functoriality}

In this and the following section we show that the results of \S\ref{crystalline} are
functorial in rational representations of the Langlands dual group and therefore
generalize to arbitrary split groups.

We use the notations and assumptions of \S\ref{satake}. In
particular $G = \mathbf{G}(L)$ with $\mathbf{G}$ an $L$-split
connected reductive group over $L$, $\widetilde{\xi} \in
X^\ast(\widetilde{\mathbf{T}})$ is the highest weight of a
$\mathbb{Q}_p$-rational representation $\rho$ of $G$ in a $K$-vector
space and $\mathbf{T}'$ is the $K$-torus dual to $\mathbf{T}$. We
assume throughout this section that $q^{1/2} \in K$. By the
normalized Satake isomorphism, any point $\zeta \in \mathbf{T}(K)$
gives rise to a $K$-valued character of the Satake-Hecke algebra
$\mathcal{H}(G,\rho_U)$ which we view as a one dimensional module
$K_\zeta$ for $\mathcal{H}(G,\rho_U)$ as in \S\ref{crystalline}. By
specialization we may form the locally algebraic $G$-representation:
\begin{equation*}
    H_{\xi,\zeta} := K_\zeta \otimes_{\mathcal{H}(G,\rho_U)}
    {\rm ind}_U^G(\rho_U)
\end{equation*}
which is of finite length and has a unique irreducible quotient
$V_{\xi,\zeta}$.

Let us first look at the case $\xi = 1$. Then the $G$-representation
$V_{1,\zeta}$ is smooth. Let $\mathbf{G}'$ be the connected
Langlands dual group over $K$ of $\mathbf{G}$. It contains
$\mathbf{T}'$ as a maximal $K$-split torus. Hence our point $\zeta
\in \mathbf{T}'(K)$ defines a $K$-split semisimple conjugacy class
in $\mathbf{G}'(K)$ which we may also view as an isomorphism class
of unramified homomorphisms ${\rm W}(\Qpbar/L) \longrightarrow
\mathbf{G}'(K)$. In the limit over $K$ the correspondence:
\begin{equation*}
    \zeta \longmapsto V_{1,\zeta}
\end{equation*}
therefore is a manifestation of the unramified local Langlands
functoriality principle (compare \cite[Chap.II \& III]{Bor}).

Going back to the case of a general $\xi$, let us assume that the
point $\zeta$ lies in the affinoid subdomain $\mathbf{T}'_{\xi,{\rm
norm}}$ of $\mathbf{T}'$. Then the corresponding character of
$\mathcal{H}(G,\rho_U)$ extends to a (continuous) character of the
Banach algebra $\mathcal{B}(G,\rho_U)$ (\S\ref{satake}) and we may
form, using the completed tensor product, the specialization:
\begin{equation*}
    B_{\xi,\zeta} := K_\zeta\; \widehat{\otimes}_{\mathcal{B}(G,\rho_U)}
    \; B_U^G(\rho_U).
\end{equation*}
It is a unitary Banach space representation of $G$.

\begin{conj}\label{nonvan}
The Banach space $B_{\xi,\zeta}$ is non-zero.
\end{conj}

If $G={\rm GL}_{d+1}(L)$ then this conjecture is a special case of
Conjecture \ref{conjecture} (see the proof of Theorem
\ref{conjunramified}).

Following the case $G={\rm GL}_{d+1}(L)$ of \S\ref{crystalline}, we
will construct, given a pair $(\xi,\zeta)$ with $\zeta \in
\mathbf{T}'_{\xi,{\rm norm}}(K)$, a family of $p$-adic Galois
representations of $\g$ with values in $\mathbf{G}'(\overline{K})$.
One naive hope would be that this family parametrizes the
topologically irreducible quotients of $B_{\xi,\zeta}$ (which, if
true, could be seen as a $p$-adic extension of unramified Langlands
functoriality). Note that, for general $G$, we are obliged to use
the {\it normalized} Satake isomorphism (there is no Tate
normalization available).

It is useful to begin in a more general setting. We view $\zeta \in
\mathbf{T}'(K) \subseteq \mathbf{G}'(K)$ as a point of the dual
group. Via the natural identification
$X^\ast(\widetilde{\mathbf{T}}) = X_\ast(\widetilde{\mathbf{T}}')$,
we view our highest weight $\widetilde{\xi}$ as a rational
cocharacter of the $K$-torus $\widetilde{\mathbf{T}}'$ dual to
$\widetilde{\mathbf{T}}$. Obviously, one has
$X_\ast(\widetilde{\mathbf{T}}') \subseteq
X_\ast(\widetilde{\mathbf{G}}')(K)$ where the right hand side
denotes the group of $K$-rational cocharacters of the connected
Langlands dual group $\widetilde{\mathbf{G}}'$ of
$\widetilde{\mathbf{G}}$. The latter satisfies:
\begin{equation*}
\widetilde{\mathbf{G}}' = \prod_{\sigma : L\hookrightarrow
K}\mathbf{G}'
\end{equation*}
so that we have:
\begin{equation*}
X_\ast(\widetilde{\mathbf{G}}') = \prod_{\sigma : L\hookrightarrow
K}X_\ast(\mathbf{G}').
\end{equation*}
Hence $\widetilde{\xi}$ gives rise to a family of cocharacters
$(\widetilde{\xi}_\sigma)_{\sigma}$ where $\widetilde{\xi}_\sigma
\in X_\ast(\mathbf{G}')(K)$.

We fix now more generally any pair $(\nu,b)$ where $\nu =
(\nu_\sigma)_\sigma$ with $\nu_{\sigma} \in X_\ast(\mathbf{G}')(K)$
and $b \in \mathbf{G}'(K)$. Let ${\rm REP}_K (\mathbf{G}')$ denote
the neutral Tannakian category of all $K$-rational representations
of $\mathbf{G}'$ and ${\rm FIC}_{L,K}$ the additive tensor category
of all filtered isocrystals over $L$ with coefficients in $K$. An
object of the latter is a triple $(D,\varphi, {\rm Fil}^\cdot D_L)$
consisting of a free $L_0 \otimes_{\mathbb{Q}_p} K$-module $D$ of
finite rank, a $\varphi_0$-linear automorphism $\varphi$ of $D$ and
an exhaustive, separated and decreasing filtration ${\rm Fil}^\cdot
D_L$ on $D_L = L \otimes_{L_0} D$ by $L \otimes_{\mathbb{Q}_p}
K$-submodules (see \S\ref{crystalline}). The pair $(\nu,b)$ gives
rise to an additive tensor functor:
\begin{equation*}
    I_{(\nu,b)} : {\rm REP}_K(\mathbf{G}') \longrightarrow {\rm FIC}_{L,K}
\end{equation*}
as follows (depending on the choice of an embedding $\sigma_0 \in
{\rm Hom}_{\mathbb{Q}_p}(L_0,K)$). Let $\rho' : \mathbf{G}'
\longrightarrow {\rm GL}(E')$ be a $K$-rational representation. We
put $D := L_0 \otimes_{\mathbb{Q}_p} E'$ and let:
\begin{equation*}
D = \prod_{\sigma :L_0\hookrightarrow K} D_\sigma
\end{equation*}
be the $L_0$-isotypic decomposition where $L_0$ acts on $D_\sigma$
through the embedding $\sigma$ (as in \S\ref{crystalline}).
Correspondingly we have the decomposition:
\begin{equation*}
1 \otimes \rho'(b) = \prod_{\sigma :L_0\hookrightarrow K}
\rho'(b)_\sigma.
\end{equation*}
We now define the Frobenius $\varphi_{\rho'(b)}$ on $D$ by
$\varphi_{\rho'(b)} := (\varphi_0 \otimes 1) \circ \varphi'$ where:
\begin{equation*}
    \varphi' | D_\sigma := \begin{cases} \rho'(b)_{\sigma_0} & \text{if
    $\sigma = \sigma_0$,}\\ {\rm id} & \text{otherwise}. \end{cases}
\end{equation*}
For each cocharacter $\nu_\sigma$ we may decompose $E'$ into its
weight spaces:
\begin{equation*}
    E' = \oplus_{i \in \mathbb{Z}}\, E'(\rho'\circ\nu_\sigma,i)
\end{equation*}
with respect to the cocharacter $\rho'\circ\nu_\sigma$ and define a
filtration on $E'$:
\begin{equation*}
    {\rm Fil}^i_{\rho'\circ\nu_\sigma} E' := \oplus_{j \geq i}\, E'(\rho'\circ\nu_\sigma,j).
\end{equation*}
Using the $L$-isotypic decomposition (as in \S\ref{crystalline}):
\begin{equation*}
D_L = \prod_{\sigma :L\hookrightarrow K} D_{L,\sigma}
\end{equation*}
and the composed $K$-linear isomorphisms $E' \longrightarrow D_L
\xrightarrow{\;{\rm pr}\;} D_{L,\sigma}$, we first transport ${\rm
Fil}^\cdot_{\rho'\circ\nu_\sigma} E'$ to a filtration ${\rm
Fil}^\cdot_{\rho'\circ\nu_\sigma} D_{L,\sigma}$ on $D_{L,\sigma}$
and define ${\rm Fil}^\cdot_{\rho'\circ\nu} D_L :=({\rm
Fil}^\cdot_{\rho'\circ\nu_\sigma} D_{L,\sigma})_{\sigma}$.

Our functor now is:
\begin{equation*}
    I_{(\nu,b)}(\rho',E') := (D,
    \varphi_{\rho'(b)}, {\rm Fil}^\cdot_{\rho'\circ\nu} D_L).
\end{equation*}

\begin{definit}\label{wadm}
The pair $(\nu,b)$ is called $L$-admissible if, for any $(\rho',E')$
in ${\rm REP}_K(\mathbf{G}')$, the filtered isocrystal
$I_{(\nu,b)}(\rho',E')$ is admissible.
\end{definit}

One checks that this definition is independent of the choice of the
embedding $\sigma_0$ which was used in the construction of the
functor $I_{(\nu,b)}$. Suppose that $(\nu,b)$ is $L$-admissible.
Then $I_{(\nu,b)}$ can be viewed as a functor:
\begin{equation*}
    I_{(\nu,b)} : {\rm REP}_K(\mathbf{G}') \longrightarrow {\rm FIC}_{L,K}^{\rm adm}
\end{equation*}
into the full subcategory ${\rm FIC}_{L,K}^{\rm adm}$ of admissible
filtered isocrystals which again is a tensor category. Moreover,
letting ${\rm Rep}_K^{\rm con}(\g)$ denote the category of finite
dimensional $K$-linear continuous representations of $\g$, the
inverse of the functor $D_{\rm cris}$ induces a faithful tensor
functor (\cite{CF}):
\begin{equation*}
    {\rm FIC}_{L,K}^{\rm adm} \longrightarrow
{\rm Rep}_K^{\rm con}(\g).
\end{equation*}
By composing these two functors, we obtain a faithful tensor functor
between $K$-linear neutral Tannakian categories:
\begin{equation*}
    \Gamma_{(\nu,b)} : {\rm REP}_K(\mathbf{G}') \longrightarrow
{\rm Rep}_K^{\rm con}(\g)
\end{equation*}
(but which is not compatible with the obvious fiber functors).
Nevertheless, as explained in \cite[\S6]{ST}, by the Tannakian
formalism this latter functor gives rise to a continuous
homomorphism of groups:
\begin{equation*}
    \gamma_{\nu,b} : \g \longrightarrow
\mathbf{G}'(\overline{K})
\end{equation*}
which is determined by the pair $(\nu,b)$ up to conjugation in
$\mathbf{G}'(\overline{K})$. So we see that any $L$-admissible pair
$(\nu,b)$ determines an isomorphism class of ``Galois parameters''
$\gamma_{\nu,b}$.

We assume in this section that $\eta$ is integral, i.e. lies in
$X^\ast(\mathbf{T})$. Via the obvious diagonal embedding:
\begin{equation*}
    X^\ast(\mathbf{T}) = X_\ast(\mathbf{T}') \hookrightarrow
    X_\ast(\widetilde{\mathbf{T}}')
\end{equation*}
we may form the product cocharacter $\widetilde{\xi}\eta =
(\widetilde{\xi}_\sigma \eta)_\sigma \in
X_\ast(\widetilde{\mathbf{G}}')(K)$. We have
$\v\circ\big((\widetilde{\xi}\eta)|_T\big) = \xi_L + \eta_L$.

\begin{thm}\label{integralmain}
Suppose that $\eta$ is integral, let $\widetilde{\xi} \in
X^\ast(\widetilde{\mathbf{T}})$ be dominant, and let $\zeta \in
\mathbf{T}'(K)$. Then there exists an $L$-admissible pair $(\nu,b)$
(and hence a Galois parameter $\gamma_{\nu,b}$) such that $\nu$ lies
in the $\widetilde{\mathbf{G}}'(K)$-orbit of $\widetilde{\xi}\eta$
and $b$ has semisimple part $\zeta$ if and only if $\zeta \in
\mathbf{T}'_{\xi,{\rm norm}}(K)$.
\end{thm}
\begin{proof}
This is a straightforward generalization of the proof of
\cite[Prop.6.1]{ST}. It might be more in the spirit of that proof
and hence helpful to note that the category ${\rm FIC}_{L,K}$ can
equivalently be described as follows. For any natural number $\ell >
0$, let ${\rm F}^\ell{\rm IC}_K$ denote the category of $\ell$-fold
filtered $K$-isocrystals whose objects are finite dimensional
$K$-vector spaces equipped with a $K$-linear automorphism and a
family of $\ell$ exhaustive, separated and decreasing filtrations by
$K$-subspaces. Fixing again an embedding $\sigma_0 \in {\rm
Hom}_{\mathbb{Q}_p}(L_0,K)$, a construction as in the definition of
$I_{(\nu,b)}$ establishes an equivalence of categories ${\rm
FIC}_{L,K} \simeq {\rm F}^{[L:\mathbb{Q}_p]}{\rm IC}_K$. We leave it
to the reader as an exercise to work out the weak admissibility
conditions for objects in the category ${\rm F}^\ell{\rm IC}_K$
(compare Proposition \ref{admissible} and its proof).
\end{proof}

We remark that the element $b$ in the statement of Theorem
\ref{integralmain} can be chosen to be a regular element in
$\mathbf{G}'(K)$ whose semisimple part is $\zeta$. Note that to make
this theorem compatible with \S\ref{crystalline}, one has to
consider $L$-admissible pairs $(\nu,b)$ such that $\nu$ lies in the
$\widetilde{\mathbf{G}}'(K)$-orbit of $(\widetilde{\xi}\eta)^{-1}$
and $b$ has semisimple part $\zeta^{-1}$.

\section{Towards a $p$-adic unramified functoriality II}

We keep the notations of \S\ref{functoriality}. We would like here
to drop the integrality assumption on $\eta$ in Theorem
\ref{integralmain}. For this, we have to introduce the additive
tensor category ${\rm FIC}_{L,K,2}$ of 2-filtered isocrystals over
$L$ with coefficients in $K$. These are triples $(D,\varphi,{\rm
Fil}^\cdot D_L)$ defined exactly as before except that the
filtration ${\rm Fil}^\cdot D_L$ is allowed to be indexed by
$\frac{1}{2}\mathbb{Z}$. The notion of weak admissibility also
extends, with literally the same definition, to the objects in this
larger category leading to the full subcategory ${\rm
FIC}_{L,K,2}^{\rm adm}$.

\begin{prop}\label{tann}
The category ${\rm FIC}_{L,K,2}^{\rm adm}$ is $K$-linear neutral
Tannakian.
\end{prop}
\begin{proof}
This is essentially contained in \cite{Tot}.
\end{proof}

Let $\mathbb{D}$ denote the protorus with character group
$\mathbb{Q}$. The elements $\nu \in (X_\ast(\mathbf{G}') \otimes
\frac{1}{2}\mathbb{Z})(K)$ can be viewed as $K$-rational
homomorphisms $\nu : \mathbb{D} \longrightarrow \mathbf{G}'$ whose
weights in any given $K$-rational representation of $\mathbf{G}'$
lie in $\frac{1}{2}\mathbb{Z}$. Hence our earlier construction of an
additive tensor functor $I_{(\nu,b)}$ makes sense for any pair
$(\nu,b) \in (X_\ast(\widetilde{\mathbf{G}}') \otimes
\frac{1}{2}\mathbb{Z})(K) \times \mathbf{G}'(K)$ producing a
functor:
\begin{equation*}
I_{(\nu,b)} : {\rm REP}_K(\mathbf{G}') \longrightarrow {\rm
FIC}_{L,K,2}.
\end{equation*}
We continue to call the pair $(\nu,b)$ $L$-admissible if this
functor has values in ${\rm FIC}_{L,K,2}^{\rm adm}$. Since $\eta \in
X^\ast(\mathbf{T})\otimes \frac{1}{2}\mathbb{Z}$, the following
variant of Theorem \ref{integralmain} holds true with literally the
same proof.

\begin{thm}\label{main}
Let $\widetilde{\xi} \in X^\ast(\widetilde{\mathbf{T}})$ be
dominant, and let $\zeta \in \mathbf{T}'(K)$. Then there exists an
$L$-admissible pair $(\nu,b)$ such that $\nu$ lies in the
$\widetilde{\mathbf{G}}'(K)$-orbit of $\widetilde{\xi}\eta$ and $b$
has semisimple part $\zeta$ if and only if $\zeta \in
\mathbf{T}'_{\xi,{\rm norm}}(K)$.
\end{thm}

For the rest of this section, we fix a pair
$(\widetilde{\xi},\zeta)$ such that $\zeta \in \mathbf{T}'_{\xi,{\rm
norm}}(K)$ and a pair $(\nu,b)$ as in Theorem \ref{main}. We have
the faithful tensor functor:
\begin{equation*}
I_{(\nu,b)} : {\rm REP}_K(\mathbf{G}') \longrightarrow {\rm
FIC}_{L,K,2}^{\rm adm}.
\end{equation*}
Our goal is to associate with $(\nu,b)$ similarly as before an
isomorphism class of Galois parameters. It is not clear how to
relate the category ${\rm FIC}_{L,K,2}^{\rm adm}$ to the category
${\rm Rep}_K^{\rm con}(\g)$, but the functor $I_{(\nu,b)}$ has
values in a particular subcategory of ${\rm FIC}_{L,K,2}^{\rm adm}$
which will turn out to be related to the category of continuous
representations of a group close to $\g$.

For any embedding $\sigma_0 \in {\rm Hom}_{\mathbb{Q}_p}(L_0,K)$,
recall that we have the $\sigma_0$-tautological fiber functor on
${\rm FIC}_{L,K,2}^{\rm adm}$ which sends the isocrystal $D$ to its
$L_0$-isotypic component $D_{\sigma_0}$. Let $\mathbb{G}$ (resp.
$\mathbb{G}_2$) denote the affine $K$-group scheme of
$\otimes$-automorphisms of the $\sigma_0$-tautological fiber functor
on ${\rm FIC}_{L,K}^{\rm adm}$ (resp. on ${\rm FIC}_{L,K,2}^{\rm
adm}$). The inclusion of Tannakian categories ${\rm FIC}_{L,K}^{\rm
adm} \subseteq {\rm FIC}_{L,K,2}^{\rm adm}$ corresponds to a
faithfully flat $K$-rational homomorphism $\mathbb{G}_2
\longrightarrow \mathbb{G}$ (\cite[Prop.2.21(a)]{DM}). Under the
functor $I_{(\nu,b)}$ (whose construction involved a choice of
$\sigma_0$), the tautological fiber functor on ${\rm
REP}_K(\mathbf{G}')$ corresponds to the $\sigma_0$-tautological
fiber functor on ${\rm FIC}_{L,K,2}^{\rm adm}$. Hence $I_{(\nu,b)}$
is induced by a homomorphism of $K$-group schemes:
\begin{equation*}
i_{(\nu,b)} : \mathbb{G} \longrightarrow \mathbf{G}'.
\end{equation*}

The homomorphism:
\begin{align*}
    L^\times & \longrightarrow K^\times\\
    a\ & \longmapsto  |a|_L\cdot N_{L/\Qp}(a)
\end{align*}
extends to the profinite completion of $L^\times$ and hence,
composed with the reciprocity map of local class field theory,
defines a continuous character:
\begin{equation*}
    \varepsilon : \g \longrightarrow K^\times
\end{equation*}
which is nothing else than the $p$-adic cyclotomic character
restricted to $\g$. Define the following object $\underline{K}$ in
${\rm FIC}_{L,K}^{\rm adm}$: its underlying module is $L_0
\otimes_{\mathbb{Q}_p} K$, its filtration is:
$$L\otimes_{\Qp}K={\rm Fil}^{-1} (L \otimes_{\mathbb{Q}_p} K)\supseteq 0={\rm Fil}^0 (L \otimes_{\mathbb{Q}_p} K)$$
and its Frobenius on $L_0 \otimes_{\mathbb{Q}_p} K$ is
$\varphi_0\otimes p^{-1}$. Equivalently, writing $L_0
\otimes_{\mathbb{Q}_p} K=\prod_{\sigma_0}K$ and rescaling, the
Frobenius is the circular permutation between the components with
multiplication by $q^{-1}$ after one round.

The following lemma is well known:

\begin{lem}
The representation $\varepsilon$ and the filtered module
$\underline{K}$ correspond to each other under the functor $D_{\rm
cris}$.
\end{lem}

The object $\underline{K}$ in ${\rm FIC}_{L,K}^{\rm adm}$ also
corresponds to a $K$-rational character:
\begin{equation*}
    \underline{\varepsilon} : \mathbb{G} \longrightarrow
\mathbb{G}_m.
\end{equation*}
We now introduce the fiber product $\mathbb{G}_{(2)}$ of affine
$K$-group schemes:
\begin{equation*}
    \xymatrix{
  \mathbb{G}_{(2)} \ar[d]_{\underline{\varepsilon}_{2}} \ar[r]
                & \mathbb{G} \ar[d]^{\underline{\varepsilon}}  \\
  \mathbb{G}_m \ar[r]^{(.)^2}
               & \mathbb{G}_m          }
\end{equation*}
and note that the horizontal arrows are faithfully flat. In
particular ${\rm FIC}_{L,K}^{\rm adm} \simeq {\rm
REP}_K(\mathbb{G})$ is a full subcategory of ${\rm
REP}_K(\mathbb{G}_{(2)})$. The kernel of the upper horizontal arrow
is central and isomorphic to the group of order two $\mu_2$. Using
the decomposition into eigenspaces with respect to $\mu_2$, one
concludes that any object $M$ in ${\rm REP}_K(\mathbb{G}_{(2)})$
decomposes uniquely as:
\begin{equation*}
    M = M_0 \oplus (\underline{\varepsilon}_{2} \otimes M_1)
\end{equation*}
with objects $M_0, M_1$ in ${\rm REP}_K(\mathbb{G})$.

In ${\rm FIC}_{L,K,2}^{\rm adm}$, the object $\underline{K}$ has a
tensor square root $\underline{K}_{2}$ defined similarly but with
$q^{-1/2}$ instead of $q^{-1}$ and with the filtration
$L\otimes_{\Qp}K={\rm Fil}^{-1/2} (L \otimes_{\mathbb{Q}_p}
K)\supseteq 0={\rm Fil}^0 (L \otimes_{\mathbb{Q}_p} K)$. We
therefore have a unique $K$-rational homomorphism:
\begin{equation*}
    \mathbb{G}_2 \longrightarrow \mathbb{G}_{(2)}
\end{equation*}
such that the composite $\mathbb{G}_2 \longrightarrow
\mathbb{G}_{(2)} \xrightarrow{\,\underline{\varepsilon}_2\,}
\mathbb{G}_m$ classifies $\underline{K}_{2}$ and the composite
$\mathbb{G}_2 \longrightarrow \mathbb{G}_{(2)} \longrightarrow
\mathbb{G}$ is the natural map. Let:
$${\rm FIC}_{L,K,(2)}^{\rm adm}$$
denote the Tannakian subcategory of ${\rm FIC}_{L,K,2}^{\rm adm}$
generated by ${\rm FIC}_{L,K}^{\rm adm}$ and the object
$\underline{K}_{2}$. By direct inspection, one checks that any
object $D$ in ${\rm FIC}_{L,K,(2)}^{\rm adm}$ decomposes uniquely
as:
\begin{equation*}
    D = D_0 \oplus (\underline{K}_{2} \otimes D_1)
\end{equation*}
with objects $D_0, D_1$ in ${\rm FIC}_{L,K}^{\rm adm}$. It follows
that the homomorphism $\mathbb{G}_2 \longrightarrow
\mathbb{G}_{(2)}$ is faithfully flat and induces an equivalence of
categories ${\rm REP}_K(\mathbb{G}_{(2)}) \simeq {\rm
FIC}_{L,K,(2)}^{\rm adm}$ so that $\underline{\varepsilon}_{2}$
corresponds to $\underline{K}_{2}$.

\begin{prop}\label{factor}
The homomorphism $i_{(\nu,b)} : \mathbb{G} \longrightarrow
\mathbf{G}'$ factorizes through $\mathbb{G}_{(2)}$.
\end{prop}
\begin{proof}
We have to show that $I_{(\nu,b)}((\rho',E'))$ lies in ${\rm
FIC}_{L,K,(2)}^{\rm adm}$ for any $(\rho',E')$ in ${\rm
REP}_K(\mathbf{G}')$.

\textit{Step 1:} We first consider the special case where the
derived group of $\mathbf{G}$ is simply connected. Then there exists
an element $\chi \in X^\ast(\mathbf{G})\otimes
\frac{1}{2}\mathbb{Z}$ such that $\chi\eta \in X^\ast(\mathbf{T})$
is integral (\cite[\S8]{Gro}). We claim that, with $(\nu,b)$, also the
pair $(\chi\nu,\chi^2(q^{1/2})b)$ is $L$-admissible. Note that
$\chi|\mathbf{T}$ viewed as a $\mathbf{T}'$-valued cocharacter in
fact has values in the connected center $\mathbf{Z}' \subseteq
\mathbf{T}' \subseteq \mathbf{G}'$ of $\mathbf{G}'$. Hence $\chi\eta
\in X_\ast(\widetilde{\mathbf{G}}')(K)$ via the diagonal embedding
and $\chi\nu$ lying in the $\widetilde{\mathbf{G}}'(K)$-orbit of
$\widetilde{\xi}\chi\eta \in X_\ast(\widetilde{\mathbf{G}}')(K)$ is
integral as well. Let now $\rho' : \mathbf{G}' \longrightarrow {\rm
GL}(E')$ be any $K$-rational representation. By additivity we may
assume that $\rho'$ is irreducible. It then follows from Schur's
lemma (\cite[Prop.\ II.2.8]{Jan}) that the image of $\rho'\circ\chi$
lies in the center of ${\rm GL}(E')$. Hence there is an $n_{E'} \in
\frac{1}{2}\mathbb{Z} \subseteq X^\ast(\mathbb{D})$ such that:
\begin{align*}
    \rho'\circ\chi : \mathbb{D} & \longrightarrow  \mathbb{G}_m = \text{center
    of ${\rm GL}(E')$} \\
    a & \longmapsto  a^{n_{E'}}
\end{align*}
and:
\begin{align*}
    \rho'\circ\chi^2 : \mathbb{D} & \longrightarrow  \mathbb{G}_m  \\
    a & \longmapsto  a^{2n_{E'}}.
\end{align*}
We conclude that:
\begin{equation*}
{\rm Fil}^\bullet_{\rho'\circ (\chi\nu)} (L \otimes_{\Qp} E') = {\rm
Fil}^{\bullet - n_{E'}}_{\rho'\circ\nu} (L \otimes_{\Qp} E')
\end{equation*}
and that $\rho'(\chi^2(q^{1/2})b) = q^{n_{E'}}\rho'(b)$. This
obviously means that:
\begin{equation*}
    I_{(\chi\nu,\chi^2(q^{1/2})b)}(\rho',E')
    = \underline{K}_{2}^{\otimes -2n_{E'}} \otimes
    I_{(\nu,b)}(\rho',E').
\end{equation*}
We in particular obtain that $I_{(\nu,b)}(\rho',E')$ lies in ${\rm
FIC}_{L,K,(2)}^{\rm adm}$.

\textit{Step 2:} For a general $\mathbf{G}$, we choose some
$L$-split $z$-extension $f : \mathbf{H} \longrightarrow \mathbf{G}$
of $\mathbf{G}$ (\cite[Lem.1.1]{Kot}). This is a surjective
$L$-rational homomorphism of $L$-split connected reductive groups
whose kernel is an $L$-split torus $\mathbf{S}$ which is central in
$\mathbf{H}$ and such that the derived group of $\mathbf{H}$ is
simply connected. By functoriality (\cite[\S I.2.5]{Bor}) we obtain
a short exact sequence of Langlands dual groups:
\begin{equation*}
    1 \longrightarrow \mathbf{G}' \xrightarrow{\ f'\ } \mathbf{H}'
    \longrightarrow \mathbf{S}' \longrightarrow 1.
\end{equation*}
Clearly we have the commutative diagram of functors:
\begin{equation*}
    \xymatrix@R=0.5cm{
    {\rm REP}_K(\mathbf{H}')  \ar[dd]_{\rm res} \ar[drr]^{I_{(f'(\nu),f'(b))}}             \\
                & & {\rm FIC}_{L,K,2}.         \\
    {\rm REP}_K(\mathbf{G}')  \ar[urr]_{I_{(\nu,b)}}              }
\end{equation*}
This implies that, with $(\nu,b)$, also $(f'(\nu),f'(b))$ is
$L$-admissible and that, in fact, we have the commutative diagram:
\begin{equation*}
    \xymatrix@R=0.5cm{
    {\rm REP}_K(\mathbf{H}')  \ar[dd]_{\rm res} \ar[drr]^{I_{(f'(\nu),f'(b))}}             \\
               & & {\rm FIC}_{L,K,2}^{\rm adm}.          \\
    {\rm REP}_K(\mathbf{G}')  \ar[urr]_{I_{(\nu,b)}}              }
\end{equation*}
By Step 1, the upper oblique arrow has values in ${\rm
FIC}_{L,K,(2)}^{\rm adm}$. But by the theory of dominant weights the
perpendicular restriction functor is surjective on objects. Hence
the lower oblique arrow has values in ${\rm FIC}_{L,K,(2)}^{\rm
adm}$ as well.
\end{proof}

On the Galois side, we may imitate the construction of
$\mathbb{G}_{(2)}$ as follows.

First, we need the {\it assumption} that any element of $\Qp^\times$
is a square in $K^\times$. We then introduce the fiber product:
\begin{equation*}
    \xymatrix{
  \g_{(2)} \ar[d]_{\varepsilon_{2}} \ar[r]
                & \g \ar[d]^{\varepsilon}  \\
  K^\times \ar[r]^{(.)^2}
                & K^\times.           }
\end{equation*}
This produces a central extension $({\rm Ext}_L)$ of the form:
\begin{equation*}
    1 \longrightarrow \{\pm 1\} \longrightarrow
    \g_{(2)} \longrightarrow
    \g \longrightarrow 1
\end{equation*}
and we have $\varepsilon_{2}^2 = \varepsilon$.

\begin{lem}
The above extension is split if and only if $[L:\Qp]$ is even.
\end{lem}
\begin{proof}
Extensions of the above form are classified by the Galois cohomology
group $H^2(\g, \mathbb{Z}/2\mathbb{Z}) =H^2(\g, \mu_2)$ which, by
local class field theory, is isomorphic to the $2$-torsion subgroup
${\rm Br}(L)_2$ in the Brauer group of $L$ and hence has order two.
One easily checks that the specific extension $({\rm Ext}_{\mathbb{Q}_p})$
is non-split and that $({\rm Ext}_L)$ is the restriction to $\g$ of
$({\rm Ext}_{\mathbb{Q}_p})$. But, again by local class field theory, the
restriction map ${\rm Br}(\mathbb{Q}_p)_2 \rightarrow {\rm Br}(L)_2$
is the multiplication by $[L:\mathbb{Q}_p]$ and hence is the zero
map if and only if $[L:\mathbb{Q}_p]$ is even.
\end{proof}

Any representation $V$ in ${\rm Rep}_K^{\rm con}(\g_{(2)})$
decomposes into its eigenspaces $V = V_+ \oplus V_-$ with respect to
the action of the subgroup $\{\pm 1\}$. Moreover, $V_+$ and
$\varepsilon_{2} \otimes V_-$ lie in ${\rm Rep}_K^{\rm con}(\g)$.
The inverse of $D_{\rm cris}$ therefore extends to a faithful tensor
functor:
\begin{equation*}
    D_{\rm cris}^{-1} : {\rm FIC}_{L,K,(2)}^{\rm adm} \longrightarrow
    {\rm Rep}_K^{\rm con}(\g_{(2)})
\end{equation*}
uniquely characterized by the property that it sends
$\underline{K}_{2}$ to $\varepsilon_{2}$. To give a more intrinsic
description, we pick a generator $t$ of $\mathbb{Z}_p(1)$ inside
$B_{\rm cris} \subseteq B_{\rm dR}$ and put:
\begin{equation*}
    B_{\rm cris,2} := B_{\rm cris}[X]/(X^2 - t) \subseteq B_{\rm dR,2} :=
    B_{\rm dR}[X]/(X^2 - t).
\end{equation*}
Let $t^{1/2}$ denote the image of $X$ in $B_{\rm cris,2}$. The
Galois action on $B_{\rm cris} \subseteq B_{\rm dR}$ extends to an
action of $\g_{(2)}$ on $B_{\rm cris,2} \otimes_{\mathbb{Q}_p} K
\subseteq B_{\rm dR,2} \otimes_{\mathbb{Q}_p} K$ by the requirement
that $\g_{(2)}$ acts on $t^{1/2}$ through the character
$\varepsilon_2$, i.e. $g(t^{1/2} \otimes 1) := t^{1/2} \otimes
\varepsilon_2(g)$ for any $g \in \g_{(2)}$. The Frobenius
$\varphi_0$ on $B_{\rm cris}$ is extended to $B_{\rm cris,2}
\otimes_{\mathbb{Q}_p} K$ by:
\begin{equation*}
\varphi_0(t^{1/2} \otimes 1) := t^{1/2} \otimes p^{1/2}.
\end{equation*}
The filtration ${\rm Fil}^\cdot B_{\rm dR}$ is extended to a
filtration on $B_{\rm dR,2}$ indexed by $\frac{1}{2}\mathbb{Z}$ by
the requirement that $t^{1/2} \in {\rm Fil}^{1/2} B_{\rm dR,2}$. For
any $V$ in ${\rm Rep}_K^{\rm con}(\g_{(2)})$ we define:
\begin{equation*}
D_{\rm cris}(V) := (B_{\rm cris,2}
\otimes_{\mathbb{Q}_p}V)^{\g_{(2)}} = \big((B_{\rm cris,2}
\otimes_{\mathbb{Q}_p}K) \otimes_K V\big)^{\g_{(2)}}
\end{equation*}
and:
\begin{equation*}
D_{\rm dR}(V) := (B_{\rm dR,2} \otimes_{\mathbb{Q}_p}V)^{\g_{(2)}}
\end{equation*}
as the invariants of the respective diagonal action of $\g_{(2)}$.
The former  is a free $L_0 \otimes_{\mathbb{Q}_p} K$-module of
finite rank equipped with a $\varphi_0$-linear Frobenius
automorphism. The latter is an $L \otimes_{\mathbb{Q}_p} K$-module
equipped with a filtration. Via the inclusion $L \otimes_{L_0}
D_{\rm cris}(V) \hookrightarrow D_{\rm dR}(V)$ this filtration
induces a filtration on $D_{\rm cris}(V)_L$ indexed by
$\frac{1}{2}\mathbb{Z}$. Hence we have a functor:
\begin{equation*}
D_{\rm cris} : {\rm Rep}_K^{\rm con}(\g_{(2)})\longrightarrow {\rm
FIC}_{L,K}.
\end{equation*}

\begin{definit}\label{cris}
A representation $V$ in ${\rm Rep}_K^{\rm con}(\g_{(2)})$ is called
crystalline if the $L \otimes_{L_0} K$-rank of $D_{\rm cris}(V)$ is
equal to ${\rm dim}_K V$.
\end{definit}

Let ${\rm Rep}_K^{\rm cris}(\g_{(2)})$ denote the full subcategory
of all crystalline representations of $\g_{(2)}$. For general $V$ we
compute:
\begin{align*}
D_{\rm cris}(V) & = D_{\rm cris}(V_+ \oplus V_-) = [(B_{\rm cris} \oplus t^{1/2}B_{\rm cris}) \otimes_{\mathbb{Q}_p} (V_+ \oplus V_-)]^{\g_{(2)}} \\
& = [(B_{\rm cris} \otimes_{\mathbb{Q}_p} V_+) \oplus (t^{1/2}B_{\rm cris} \otimes_{\mathbb{Q}_p}V_-)]^{\g} \\
& = D_{\rm cris}(V_+) \oplus (\underline{K}_{2}^\ast \otimes D_{\rm
cris}(\varepsilon_2 \otimes V_-)).
\end{align*}
It follows that $V$ is crystalline if and only if $V_+$ and
$\varepsilon_2 \otimes V_-$ are crystalline in the usual sense. As
an immediate consequence of the main result of \cite{CF} we
therefore obtain the following result.

\begin{prop}\label{equiv}
The functor $D_{\rm cris}$ restricts to an equivalence of tensor
categories:
\begin{equation*}
{\rm Rep}_K^{\rm cris}(\g_{(2)}) \xrightarrow{\ \sim\ } {\rm
FIC}_{L,K,(2)}^{\rm adm}.
\end{equation*}
\end{prop}

By Proposition \ref{factor}, the composite:
\begin{equation*}
    \Gamma_{(\nu,b)} := D_{\rm cris}^{-1} \circ I_{(\nu,b)} :
    {\rm REP}_K(\mathbf{G}') \longrightarrow {\rm Rep}_K^{\rm con}(\g_{(2)})
\end{equation*}
is well defined as a faithful tensor functor between $K$-linear
neutral Tannakian categories. It gives rise, as before, to an
isomorphism class of ``Galois parameters'':
\begin{equation*}
    \gamma_{\nu,b} : \g_{(2)} \longrightarrow
\mathbf{G}'(\overline{K})
\end{equation*}
which is uniquely determined by the $L$-admissible pair $(\nu,b)$.

\end{document}